\documentclass[11pt,reqno]{amsart}
\usepackage{amsmath}
\usepackage{txfonts}
\usepackage{amssymb}
 \usepackage{xcolor}
\usepackage{amsmath,amsthm}
\usepackage{CJK,amsfonts}
\usepackage[body={7.5in,8.75in}, top=1.2in, right=0.9in,left=0.9in]{geometry}
\usepackage{hyperref}
\allowdisplaybreaks[4] \numberwithin{equation}{section}

\newtheorem{thm}{Theorem}
\newtheorem{theorem}[thm]{Theorem}

\newtheorem{lemma}[thm]{Lemma}
\newtheorem{remark}[thm]{Remark}
\newtheorem{proposition}[thm]{Proposition}

\def\squarebox#1{\hbox to #1{\hfill\vbox to #1{\vfill}}}

\newcommand{\be}{\begin{equation}}
\newcommand{\ee}{\end{equation}}
\newcommand{\bea}{\begin{eqnarray}}
\newcommand{\eea}{\end{eqnarray}}
\newcommand{\bd}{\begin{displaymath}}
\newcommand{\ed}{\end{displaymath}}
\begin{document}
\title[On certain properties of perturbed Freud-type  weight: a revisit]{On certain properties of perturbed Freud-type  weight: a revisit}
\author{  Abey S. Kelil${}^1$, Appanah R. Appadu${}^1$ and Sama Arjika${}^{2}$}
\dedicatory{\textsc}
\thanks{ ${}^{1}$Nelson Mandela University.
${}^2$Department of Mathematics and Informatics, University of Agadez, Niger.  }
\thanks{Email:  rjksama2008@gmail.com}



\begin{abstract}
In this paper, monic polynomials orthogonal with  deformation of the Freud-type weight function 	are considered. These polynomials 	fullfill linear differential equation with some polynomial coefficients in their holonomic form. The aim of this work is explore certain characterizing properties of perturbed Freud type polynomials such as  nonlinear recursion relations, finite moments,  differential-recurrence and differential relations satisfied by the recurrence coefficients as well as the corresponding semiclassical orthogonal polynomials. We note that  the obtained differential equation fulfilled by the considered semiclassical polynomials are used to study an electrostatic interpretation for the
	distribution of zeros based on the original ideas of Stieltjes.
\end{abstract}

\maketitle

\section{Introduction}
  Suppose we have a family of polynomials $\{\psi_m(x)\}_{m=1}^{\infty}$ which are monic of degree $m$ and that are orthogonal with respect to the positive weight $w(x)$ on the interval $[c,\,d]$, i.e.,
$$\langle\psi_m, \psi_k \rangle_w=\int_c^d\psi_m(x)\psi_k(x)w(x){\rm d}x=\Gamma_m\delta_{m,k},\quad m,k=0,1,2,\cdots,$$
where $\Gamma_m>0$ denotes the normalization constant \cite{refChihara,refSzego}.
This value can be  obtained from the square of the weighted    $L^2$-norm of $\psi_m(x)$ over $[c,d]$. Monic polynomial representation takes the form $$\psi_n(x)=x^n+p(n)x^{n-1}+\ldots$$
It is known that
$\det\left(x_j^{i-1}\right)_{i,j=1}^{N}=\prod _{1\leq i<j\leq N}(x_i-x_j)=\det\left(\psi_{i-1}(x_j)\right)_{i,j=1}^N .$ The polynomials $\psi_n(x)$ can be generated by the Gram-Schmidt orthogonalization process \cite{refChihara,refIsmaila}.

\noindent As it is known in \cite{refChihara,refIsmaila,refSzego}, classical orthogonal polynomials obey  Pearson's differential equation
\begin{equation}\label{eq:Pearson}
\dfrac{\,d\left(\lambda(x)w(x)\right)}{\mathrm{d}x}
=\tau(x)w(x),\end{equation}
where the polynomials $\lambda(x)$ and $\tau(x)$ are of  degrees two and one respectively.
Whereas polynomials for which the weight fullfills Eq.  \eqref{eq:Pearson} with deg$(\lambda)\geq 2$ or deg$(\tau)\neq 1$
are said to be Semi-classical orthogonal polynomials \cite{refhendriksen1985semi}.

\noindent For deformed orthogonality weight, if the moments exist and the corresponding monic orthogonal polynomials $\psi_n(z)$ for $n=0,1,2,\ldots$ obey linear recursive relation
\begin{align*}
\begin{cases}
&
z\psi_n(z)=\psi_{n+1}(z)+\gamma_n\psi_{n-1}(z)+\alpha_n\psi_n(z),
\\&\psi_0(z)=1,\quad\gamma_0 \psi_{-1}(z)=0.
\end{cases}
\end{align*}
\noindent The following relations  in \cite{ChenIsmail1997}  are valid for a semiclassical weight $w$ with  $w(a)=w(b)=0$.
\begin{lemma}{\cite{ChenIsmail1997}.}
	Suppose that  $v(x)=-\ln w(x)$ has a derivative in some Lipschitz order with a positive exponent \cite{refSearcoid}.
	The differential-difference coefficients obey the following formulas:
	\begin{align}
	\psi'_n(z) &=\gamma_n \mathcal{A}_n(z)\psi_{n-1}(z)-\mathcal{B}_n(z)  \psi_n(z), \label{UP}\\
	\psi'_{n-1}(z) &=-\mathcal{A}_{n-1}(z)\psi_n(z)+\left[\mathcal{B}_n(z)+v'(z)\right]\psi_{n-1}(z),\label{DOWN}
	\end{align}
	where
	\begin{align}
	&\mathcal{A}_n(z):=\frac{1}{\Gamma_n}\int_a^b\frac{v^\prime(z)-v^\prime(\tau)}{z-\tau}\psi^2_n(y)w(\tau){\rm d}\tau,\label{A}\\
	&\mathcal{B}_n(z):=\frac{1}{\Gamma_{n-1}}\int_a^b\frac{v^\prime(z)-v^\prime(\tau)}{z-\tau}\psi_n(\tau)\psi_{n-1}(\tau)w(\tau){\rm d}\tau.\label{B}
	\end{align}
\end{lemma}

\begin{lemma}{\cite{ChenIsmail1997}.}
	The coefficients $\mathcal{A}_n(z)$ and $\mathcal{B}_n(z)$ defined by Eq. \eqref{A} and Eq. \eqref{B} obeys 
	\begin{subequations}\label{compAb}
		\begin{align}
		\begin{cases}
		&\mathcal{B}_{n+1}(z)+\mathcal{B}_n(z)=- v'(z)+(z-\alpha_n)\mathcal{A}_n(z),~~~~~~ ~~~~~~~\qquad\qquad
		\quad
		~~~~~~~~~~~~~~~~~~~(M_1)\notag\\
		&1+  (z-\alpha_n)[\mathcal{B}_{n+1}(z)-\mathcal{B}_n(z)]=-\gamma_n\mathcal{A}_{n-1}(z)+\gamma_{n+1}\mathcal{A}_{n+1}(z).
		~~~~~\qquad\qquad~~~~(M_2)\notag
		\end{cases}
		\end{align}
	\end{subequations}
\end{lemma}
\noindent We also mention another supplementary condition, that involves $\sum _{j=0}^{n-1}\mathcal{A}_{j}(z)$ and we will denote it by  $(M_2')$ as this relation helps to obtain recurrence coefficients $\alpha_n$ and $\gamma_n$, as
\begin{equation}\label{}
v'(z)\mathcal{B}_n(z)+\sum_{j=0}^{n-1}\mathcal{A}_{j}(z)+\mathcal{B}^2_n(z)=\gamma_n\mathcal{A}_n(z)\mathcal{A}_{n-1}(z).~~~~~~~~~~~~~~~~~~~~~~~~~~(M_2')\notag
\end{equation}
Eq. $(M_2')$ can be perceived as an equation for  $\sum_{j=0}^{n-1}\mathcal{A}_{j}(z)$. See, for instance, \cite{BC2009,ChenIts2010}.

\noindent
The  differential equation fullfilled by $\psi_n(z)$ is generated by eliminating $\psi_{n-1}(z)$ from ladder  operators, and it is given as
\begin{equation}\label{H}
\psi_n''(z)-\left( v'(z)+\frac{\mathcal{A}'_n(z)}{\mathcal{A}_n(z)}\right) \psi_n'(z)+\left(
\mathcal{B}'_n(z)-\mathcal{B}_n(z)\frac{\mathcal{A}'_n(z)}{\mathcal{A}_n(z)}+\sum_{j=0}^{n-1}\mathcal{A}_{j}(z)\right)\psi_n(z)=0,
\end{equation}
where $\sum_{j=0}^{n-1}\mathcal{A}_{j}(z)$ is obtained from $(M'_2)$.
\begin{lemma}
	Suppose we have a symmetric semi-classical weight $W_{\sigma}(x;t)=\exp(tx^2)w_0(x)$,
	with $t\in\mathbb{R}$ such that  the moments of for $w_0$ is finite. The recursive coefficient $\gamma_{n}(t)$ obeys the Volterra, or the Langmuir lattice, equation \cite{van2017orthogonal}
	\begin{align}
	\dfrac{\,d{\gamma_n(t)}}{\,d{t}} = \gamma_{n}(t)\left(\gamma_{n+1}(t)-\gamma_{n-1}(t)\right).\label{eq:langlat}
	\end{align}
\end{lemma}
\begin{proof} {See, for example,  \cite[Theorem 2.4]{van2017orthogonal}.}
\end{proof}

\noindent
In this paper, we  consider to study semiclassical perturbed Freud-type measure
\begin{align}\label{sexticFreud}
\begin{cases}
&
\,d{\mu_{\sigma}}(x) = W_{\sigma}(x;t)~\mathrm{d}x = |x|^{2\sigma +1}\exp\left(-[ cx^6+t(x^4 -x^2)]\right)~\mathrm{d}x ,
\\&\sigma>0, ~c>0, ~t\in \mathbb{R},
\end{cases}
\end{align}
involving parameters $t,\, \sigma$, which will be used to represent the polynomials and in the $L^2$ norm. For simplicity, we may not sometimes display the parameters in the polynomials.  

\noindent The  motives for the choice of the  perturbed orthogonality measure in \eqref{sexticFreud} is  as follows:- First, from
some of the classical orthogonal polynomials, a new class of semiclassical (non-classical) orthogonal polynomials can be obtained by means of slight modifications on their orthogonality measure \cite{refPPNevai,refnevai1986geza}. Such measure deformation usually results in some difficulties,  most of which have not been handled yet as noted in \cite{refPPNevai,refnevai1986geza}. Motivated by the works of P. Nevai et al. \cite{refnevai1986geza}, a slight modification of a new orthogonality measure on non-compact support presents a new class of orthogonal polynomials if certain characterizing properties associated with the considered polynomials are successfully obtained. Secondly,  the choice of modified  Freud-type measure is reasonable in the sense that this orthogonality measure emanates from quadratic transformation and Chihara's symmetrization of the modified Airy-type measure (cf. \cite{refChihara} for symmetrization process).
This also leads to an investigation of certain fresh properties such as nonlinear differential-recurrence and differential equations satisfied by the recurrence coefficients as well as the perturbed polynomials themselves.  The results obtained also motivate considerable applications; for instance, in modeling nonlinear phenomena, Soliton Theory and Random matrix theory \cite{ChenIts2010} and in the crystal structure in solid-state physics, to mention a few.

\section{Semiclassical perturbed Freud-type polynomials}
\noindent
Semiclassical perturbed Freud polynomials $\lbrace\mathcal{S}_{n}(x;t)\rbrace_{n=0}^{\infty}$ on $\mathbb{R}$
are real polynomials with their orthogonality weight given by
\begin{align*}
\begin{cases}
&
\,d{\mu_{\sigma}}(x) = W_{\sigma}(x;t)~\mathrm{d}x = |x|^{2\sigma +1}\exp\left(-[ cx^6+t(x^4 -x^2)]\right)~\mathrm{d}x ,
\\&\sigma>0, ~c>0, ~t\in \mathbb{R},
\end{cases}
\end{align*}
 and  the orthogonality condition is given by
{\small\begin{align}\label{orthosextic}
	\langle \mathcal{S}_{n}, \mathcal{S}_{m} \rangle_{W_{\sigma}}=
	\int_{-\infty}^{\infty}\mathcal{S}_{n}(x;t)~\mathcal{S}_{m}(x;t)~W_{\sigma}(x;t) \,\mathrm{d}x= \hat{\Gamma}_{n}~\delta_{mn},
	\end{align}}
where $\delta_{mn}$ denotes the Kronecker delta function.
It follows from Eq. \eqref{orthosextic}  that the recursion relation takes the form
\begin{equation}
\label{symmetricrecur}
\begin{cases}
&	\mathcal{S}_{n+1}(x;t)=-\gamma_{n}(t;\sigma)~\mathcal{S}_{n-1}(x;t)+x\mathcal{S}_{n}(x;t), ~~n\in \mathbb{N},
\\& \mathcal{S}_{0}:=1~~~{\rm and}~~~\gamma_{0} \mathcal{S}_{-1}:=0.
\end{cases}
\end{equation}
If we multiply  Eq. \eqref{symmetricrecur}  with  {\small$\mathcal{S}_{n-1}(x;t)W_{\sigma}(x;t)$} and then  integrate with respect to $x$  and using  orthogonality given in Eq. \eqref{orthosextic}, we obtain
{\small\begin{align}\label{betafromularecur}
	\gamma_{n}(t;\sigma) = \frac{1}{\hat{\Gamma}_{n-1}(t)}
	\langle x\mathcal{S}_{n}, \mathcal{S}_{n-1} \rangle_{W_{\sigma}}
	=\frac{\hat{\Gamma}_{n}(t)}{\hat{\Gamma}_{n-1}}>0.
	\end{align}}
\noindent
Observe that {\small$\mathcal{S}_{n}(x;t)$} comprises the terms $x^{n-r}$, $~r\leq n$ and is symmetric so that
\begin{equation*}
\begin{cases}
&	\mathcal{S}_{n} (-x;t)=(-1)^{n} \mathcal{S}_{n}(x;t) ,
\\&\mathcal{S}_{n}(0;t) ~\mathcal{S}_{n-1}(0;t)=0,
\end{cases}
\end{equation*}
as the weight  $W_{\sigma}(x;t)$ is  even on $\mathbb{R}$.
Using monic representation of considered polynomials {\small$\mathcal{S}_{n}(x;t)$}, associated with {\small$W_{\sigma}(x;t)$}, we have that
{\small\begin{align}\label{symmetricpolyexpres}
	\mathcal{S}_{n}(x;t)=x^{n} +\chi(n;t)~x^{n-2}+\hdots + \mathcal{S}_{n}(0;t),
	\end{align}}
which can be expressed equivalently as  \cite{refChihara},
{\small\begin{equation*}
	\begin{cases}
	&	\mathcal{S}_{2i}(x;t)=x^{2i}+\chi(2i;t)~x^{2i-2}+\cdots+\mathcal{S}_{2i}(0),
	\\&
	\mathcal{S}_{2i+1}(x;t)=x^{2i+1}+\chi(2i+1;t)~x^{2i-1}+\cdots+{\rm s.}x
	=x\left(x^{2i}+\chi(2i+1;t)~x^{2i-2}+\cdots+{\rm s}\right),
	\end{cases}
	\end{equation*}}
where $s\in\mathbb{R}$.
\noindent By substituting   Eq. \eqref{symmetricpolyexpres} into  Eq. \eqref{symmetricrecur}, we obtain
\begin{align}\label{diffbetan}
\begin{cases}
&
\gamma_n(t) = \chi(n;t)-\chi(n+1;t),
\\&
\chi(0):=0.
\end{cases}
\end{align}
Imposing a telescoping iteration of terms of Eq. \eqref{diffbetan}  gives
\[\sum_{k=0}^{n-1} \gamma_k(t,\sigma)=-  \chi(n;t).\]

\section{Certain properties of the considered semiclassical polynomials}
\noindent
In this section, we explore certain characterizing properties for perturbed semi-classical Freud-type polynomials.
\subsection{Finite moments}
For certain semiclassical weights, it is known in \cite{clarkson2014relationship,clarkson2016generalized,kelil2018properties} that the moments make  link between the weight function and the theory of integrable equations, in particular, Painlev\'{e}-type equations \cite{van2017orthogonal}.

\noindent
\begin{theorem}\label{finitmom}
	Suppose $x,t\in\mathbb{R}$ and $c,\sigma>0$. The first moment $\eta_0(t;\sigma)$ associated with the weight \eqref{orthosextic} is finite.
\end{theorem}
\begin{proof} For the weight given in Eq. \eqref{sexticFreud}, the moment $\eta_0(t;\sigma)$ takes the form
	\begin{align}\label{momfnte-one}
	\eta_0(t;\sigma) = \int_{-\infty}^{\infty} W_{\sigma}(x;t) 
	~\mathrm{d} x
	= 2 \int_{0}^{\infty}
	W_{\sigma}(x;t) 
	~\mathrm{d} x .
	\end{align}
	For $\sigma >0$ and $c>0$,  the function $ W_{\sigma}(x;t) =x^{2\sigma +1}\exp\left(-[ cx^6+t(x^4 -x^2)]\right)$ is continuous on  $[0,\infty)$, and hence is integrable on $[0,\mathcal{K}]$ for any $\mathcal{K}>0$.  
	\noindent
	In order to show $ \int_{\mathcal{K}}^{\infty} W_{\sigma}(x;t)~\mathrm{d} x $ is finite, we first note that $\lim_{x\rightarrow \infty} x^{2} W_{\sigma}(x;t)=0$; that is, there exists an $N>0$ such that  $x^2 W_{\sigma}(x;t)<1$ whenever $x>N$ by definition. As  $ \int_{N}^{\infty}\frac{d{x}}{x^2}<\infty$, it follows, for $N>0$, that $\int_{N}^{\infty} W_{\sigma}(x;t)~\mathrm{d} x <\infty$,
	particularly when $N=\mathcal{K}$. Hence,  $ \int_{0}^{\infty} W_{\sigma}(x;t)~\mathrm{d} x <\infty$.
\end{proof}
\noindent The following result presents some conditions 
for  differentiation and integration order for functions of two variables \cite{jost2006postmodern}.
\begin{lemma}{\cite[Theorem 16.11]{jost2006postmodern}}.\label{jurgenjostcc}
	
	\noindent
	Let $J=(a,b)\subset\mathbb{R}$ be an open interval and $g:\mathbb{R}\times J\rightarrow\mathbb{R}$. Assume that
	\begin{enumerate}
		
		\item[(i)]
		$g(x,t)$ has a derivative on $\mathbb{R}$ with respect to $t$ for almost all $x\in\mathbb{R}$,
		\item[(ii)]
		for every fixed $t\in J$,  $\int_{-\infty}^{\infty}g(x,t)~\mathrm{d} x <\infty$,
		
		\item[(iii)]
		$\exists$ an integrable function $h:\mathbb{R} \rightarrow\mathbb{R}$ such that 
		$\forall t\in J$,
		$  \left\vert \dfrac{\partial{g(x,t)}}{{\partial{t}}} \right\vert \leq h(x), $
		which is true for almost all $x\in\mathbb{R}$.
	\end{enumerate}
	It then follows that \[	\dfrac{\,d}{\,d{t}}  \int_{-\infty}^{\infty}g(x,t)~\mathrm{d} x   =  \int_{-\infty}^{\infty}  \dfrac{\partial{g(x,t)}}{{\partial{t}}} ~ \mathrm{d} x.\]
\end{lemma}
\noindent
The following result shows how moments of high order behave for the weight function in Eq.  \eqref{sexticFreud}.
\begin{theorem}
	For $n\in \mathbb{N}_0$, the moments associated with the perturbed Freud weight given in \eqref{sexticFreud} obey 	the following formulations
	\begin{align}
	\label{mainmomentAA}
	\begin{cases}
	\eta_{2n}(t;\sigma)&=
	\ \dfrac{\,d^n}{dt^n}
	\int_{-\infty}^{\infty}|x|^{2\sigma +1}\exp\left(-[ cx^6+t(x^4 -x^2)]\right)~\mathrm{d} x  
	\\& =
	\sum_{k=0}^{n} (-1)^{n+k}\binom{n}{k}\eta_{4n-2k}(t;\sigma)= \dfrac{\,d^n}{dt^n}\eta_0(t;\sigma),
	\\
	\eta_{2n+1}(t;\sigma) & =0.
	\end{cases}
	\end{align}
	
\end{theorem}
\begin{proof}
	Taking into account  the  weight in Eq. \eqref{sexticFreud} is even on $\mathbb{R}$,  let's  take  Freud-type weight defined on the positive x-axis; that is,
	\begin{align*}
	W_{\sigma}(x;t):=x^{ 2\sigma +1}\exp\left(-[ cx^6+t(x^4 -x^2)]\right),~~x\in (0,\infty),~~\sigma>0,~~t\in J \subset\mathbb{R}.
	\end{align*}
	One can see that  $W_{\sigma}$ is a rapidly decreasing function \cite{jost2006postmodern}.

\noindent Using  Theorem \ref{finitmom}, we can easily see that 
	\begin{align}\label{keyjj}
	\dfrac{\partial{ W_{\sigma}(x;t)}}{{\partial{t}}}   =   (x^4-x^2)~x^{2\sigma +1}\exp\left(-[ cx^6+t(x^4 -x^2)]\right),
	\end{align}
	\nonumber  is continuous on $\mathbb{R}^{+}$.
	For  $t \leq 0$ and
	$x\in (1,\infty)$, we have  
	$\exp\left(t(x^4-x^2)\right)\leq 1,$ since $ty^2\leq 0$ for $y\in\mathbb{R}$. Thus,
	\begin{align}\label{bounding2b}
	\left\vert \dfrac{\partial{			W_{\sigma}(x;t)}}{{\partial{t}}}\right\rvert =  \left\vert    x^{2\sigma +1}  (x^4-x^2)\exp\left(-[ cx^6+t(x^4 -x^2)]\right) \right\vert \leq  x^{2\sigma +k}\exp\left(-cx^6\right):=G(x),
	\end{align}
	for some bounding $k\in\mathbb{R}^{+}$ and $ \sigma>0,$ with
	$$ \int_{0}^{\infty} G(x)~\mathrm{d} x = \int_{0}^{\infty} x^{2\sigma +k}\exp\left(-cx^6\right)~\mathrm{d} x   = \frac{1}{6}\left(\frac{1}{c}\right)^{\frac{\sigma+4}{k}}~\Gamma\left(\dfrac{2\sigma +8}{6}\right) <\infty,$$
	where $\Gamma(z)$ denotes the Gamma function.

\noindent It then follows from Eq. \eqref{keyjj} that
	\[
	\left\vert\dfrac{\partial{W_{\sigma}(x;t)}}{{\partial{t}}} \right\vert =  \left\vert   x^{2\sigma +3}\exp\left(-[ cx^6+t(x^4 -x^2)]\right) \right\vert \leq  x^{2\sigma +3}\exp\left(-cx^6 +Ax^2\right):=K(x),
	\]for  $t\in [0,A]$, $A\in\mathbb{R}^{+}$ and  $K(x)$ is  integrable for $x\in\mathbb{R}^{+}$.
	\noindent We see that  all the conditions of Lemma \ref{jurgenjostcc} are fulfilled so that  Eq. \eqref{mainmomentAA} can be proved using the principles of  mathematical induction. For $n=1$, we have
	\begin{align}\nonumber
	\dfrac{\,d}{\,d{t}}\eta_{0}(t,\sigma) &=  \dfrac{\,d}{\,d{t}} \int_{-\infty}^{\infty} |x|^{2\sigma +1}\exp\left(-[ cx^6+t(x^4 -x^2)]\right)~\mathrm{d} x
	\nonumber\\&=  (-1)\int_{-\infty}^{\infty}(x^4 -x^{2}) W_{\sigma}(x;t)~\mathrm{d} x
	= (-1) \left(\eta_4(t,\sigma)-\eta_2(t,\sigma)\right).
	\end{align}
	We suppose, for inductive assumption, that
	\begin{align}\label{inductiveproof1}
	\dfrac{\,d^{n}}{\,d{t^{n}}}\eta_{0}(t,\sigma) =
	\sum_{k=0}^{n} (-1)^{n+k}\binom{n}{k}\eta_{4n-2k}(t;\sigma):=\eta_{2n}(t,\sigma).
	\end{align}
	We need to show that
	\begin{align}\label{inductiveproof}
	\eta_{2n+2}(t,\sigma)= \dfrac{\,d^{n+1}}{\,d{t^{n+1}}}\eta_{0}(t,\sigma).
	\end{align}
	We note that   $x^{2n+2\sigma +1}\exp\left(-[ cx^6+t(x^4 -x^2)]\right),~x\in\mathbb{R}^{+},~ t\in J,$ also obeys the conditions of Lemma \ref{jurgenjostcc}.
	Then, by applying binomial expansion, we have
	\begin{align}\nonumber
	\dfrac{\,d^{n+1}}{\,d{t^{n+1}}}\eta_{0}(t,\sigma) &= \dfrac{\,d}{\,d{t}} \left(   \dfrac{\,d^{n}}{\,d{t^{n}}}\eta_{0}(t,\sigma) \right)
	\\ \nonumber&
	= \dfrac{\,d}{\,d{t}} \int_{\mathbb{R}} (-1)^n \left(x^4-x^2\right)^{n} W_{\sigma}(x;t)~\mathrm{d} x  
	=  \int_{\mathbb{R}} (-1)^n (-1) \left(x^4-x^2\right)^{n+1} W_{\sigma}(x;t)~\mathrm{d} x  
	\\ \nonumber&
	=  \sum_{k=0}^{n} (-1)^{n+1}\binom{n+1}{k}
	\int_{-\infty}^{\infty}  \left(x^4\right)^{n+1-k}\left(-x^2\right)^{k}~  W_{\sigma}(x;t) ~\mathrm{d} x
	\\ \nonumber&=
	\sum_{k=0}^{n} (-1)^{n+k+1}\binom{n+1}{k}
	\eta_{4n+4-2k}(t;\sigma)
	=\eta_{2n+2}(t,\sigma)\equiv  \eta_{0}(t; n+\sigma+1).
	\end{align}
	\noindent Besides,  moments of odd order vanish; i.e.,  
	\begin{equation}\label{hhkarm}
	\eta_{2n+1}(t;\sigma) = \int_{-\infty}^{\infty} x^{2n+1}\,  W_{\sigma}(x;t)\mathrm{d}x =0, ~~ n\in \mathbb{N},
	\end{equation}
	as the expression in the above integral  is an odd function.
\end{proof}

\subsection{Concise formulation}
\noindent The following result gives concise formulation for perturbed Freud-type polynomials $\mathcal{S}_n(x;t)$. For a similar result, \cite[Lemma 3.2]{refIsmailMansour} .
\begin{lemma}
	Suppose we have the perturbed  Freud-type weight given in \eqref{sexticFreud}.
	Concise formulation of the corresponding polynomials, in terms of  recurrence coefficient $\gamma_j(t;\sigma)$, is given  by
	{\small\begin{subequations}\label{coef}
			\begin{align}\label{polyforAc}
			\begin{cases}
			&			\mathcal{S}_{q}(x;t)= \displaystyle \sum_{k=0}^{\lfloor\frac q2\rfloor}\Psi_{k}(q)\medspace x^{q-2k},
			\\&
			\Psi_0(q)=1, \quad \text{ for}\quad k\in\{1,2,\dots,\lfloor\frac q2\rfloor\},~~q\in \mathbb{N},
			\end{cases}
			\end{align}
			where 
			\begin{equation}
			\label{maincolsedfora}
			\Psi_{k}(q) =(-1)^{k}  \sum_{j_{1}=1}^{q+1-2k} \gamma_{j_{1}}(t;\sigma)
			\displaystyle \sum_{j_{2}=j_{1}+2}^{q+3-2k} \gamma_{j_{2}}(t;\sigma)
			\displaystyle \sum_{j_{3}=j_{2}+2}^{q+5-2k} \gamma_{j_{3}}(t;\sigma) \cdots   \displaystyle \sum_{j_{k}=j_{k-1}+2}^{q-1} \gamma_{j_{k}}(t;\sigma).
			\end{equation}
	\end{subequations}}
\end{lemma}
\begin{proof}
	Since the perturbed Freud-type polynomials $\mathcal{S}_{q}(x;t)$ are symmetric  and monic of degree $q$, and for a fixed {\small$t\in \mathbb{R}$}, we have {\small$\mathcal{S}_q(-x)= (-1)^{q}\mathcal{S}_q(x),$} so that
	\begin{align}
	\mathcal{S}_{2q}(x;t) =  \sum_{j=0}^{q} g_{2q-2j}~x^{2q-2j};\qquad
	\mathcal{S}_{2q+1}(x;t) =  \sum_{j=0}^{q} g_{2q-2j +1}~x^{2q-2j +1},
	\end{align}
	where  {\small$g_{q-2k}=\Psi_{k}(q)$}  with {\small$\Psi_0(q)=1$} and {\small$\Psi_k(q)=0$} for {\small$k>\lfloor\frac q2\rfloor$}. If we substitute  Eq. \eqref{polyforAc} into  Eq. \eqref{symmetricrecur} and  if we compare the coefficients of $x$, we obtain
	\begin{align}\label{functionalrelaAs}
	\begin{cases}
	&		\Psi_k(q+1) -  \Psi_{k}(q) =  - \gamma_{q}(t;\sigma)   \Psi_{k-1}(q-1),
	\\&\Psi_0(q)=1.
	\end{cases}
	\end{align}
	Eq. \eqref{maincolsedfora} can be proved by employing induction on {\small$k$}. For $k=1$, we  see that
	{\small\begin{align}\label{Eqmin11}
		\Psi_1(q) -  \Psi_{1}(q-1) =  - \gamma_{q-1},
		\end{align}}
	By employing  a  telescoping sum of terms in  Eq. \eqref{Eqmin11}, we obtain
	{\small \[\Psi_1(q) =  -\sum_{j_1=0}^{q-1}\gamma_{j_{1}}(t;\sigma), \medspace \forall q\geq 1.\] }
	\noindent
	Let's assume that,  for every {\small$q\in \mathbb{N}$}, Eq. \eqref{maincolsedfora} holds true for values up to {\small$k-1$}, i.e.,
	{\small \begin{align}\label{maincolsedfoccr}
		\Psi_{k-1}(n) =(-1)^{k-1} \displaystyle \sum_{j_{1}=1}^{q+3-2k} \gamma_{j_{1}}(t;\sigma)
		\displaystyle \sum_{j_{2}=j_{1}+2}^{q+5-2k} \gamma_{j_{2}}(t;\sigma)
		\displaystyle \sum_{j_{3}=j_{2}+2}^{q+7-2k} \gamma_{j_{3}}(t;\sigma) \cdots   \displaystyle \sum_{j_{k-1}=j_{k-2}+2}^{q-1} \gamma_{j_{k-1}}(t;\sigma).  
		\end{align}}
	Eq. \eqref{functionalrelaAs} can be repeatedly used to obtain 
	{\small{\begin{align}\label{sortofrecurrclosed}
			\nonumber
			\Psi_k(q)& =   \Psi_{k}(q-1)   - \gamma_{q-1}\Psi_{k-1}(q-2),  
			\nonumber \\&
			= \Psi_{k}(q-2)   - \gamma_{q-2}\Psi_{k-1}(q-3) - \gamma_{q-1}\Psi_{k-1}(q-2),    
			\nonumber \\&
			= \Psi_{k}(q-3)   - \gamma_{q-3}\Psi_{k-1}(q-4)  - \gamma_{q-2}\Psi_{k-1}(q-3) - \gamma_{q-1}\Psi_{k-1}(q-2),    
			\nonumber \\&
			\qquad \qquad\qquad\vdots
			\nonumber \\&
			= -\gamma_{2k-1}\Psi_{k-1}(2k-2)  -\gamma_{2k}\Psi_{k-1}(2k-1)
			-\cdots  - \gamma_{q-2}\Psi_{k-1}(q-3) - \gamma_{q-1}\Psi_{k-1}(q-2).
			\end{align}}}
	\noindent Substituting Eq. \eqref{maincolsedfoccr} into  Eq. \eqref{sortofrecurrclosed} yields  Eq. \eqref{maincolsedfora} and hence the required result.
\end{proof}
\noindent Lemma 5 is alternately given as follows.
\begin{proposition}\label{niceform}
	The following formulation also holds for monic perturbed Freud-type polynomials $\mathcal{S}_{q}(x;t)$:
	{\small \begin{align}\nonumber
		\mathcal{S}_{q}(x;t)= x^{q} +  \sum_{r=1}^{\lfloor\frac q2\rfloor} (-1)^{r} \left( \displaystyle \sum_{k\in W(q,r)} \gamma_{k_1} \gamma_{k_2}\cdots \gamma_{k_{r-1}}\gamma_{k_{r}}\right) \medspace x^{q-2r},
		\end{align}}
	where
	$W(q,r) = \{k\in \mathbb{N}^{r} ~|~ k_{j+1} \geq k_j +2~~\text{for}~~1\leq j\leq r-1,~~1\leq k_1,k_r<q\},$
	and
	$	\lfloor \frac{q}{2} \rfloor
	= \begin{cases}
	&\frac{q}{2}, \quad \text{$q$ is even},
	\\ &\frac{q-1}{2}, \quad \text{$q$ is odd}.
	\end{cases}
	$
\end{proposition}

\subsection{Normalization constant }
\noindent
The normalization constant {\small$\hat{\Gamma}_{m}$} in  Eq. \eqref{orthosextic}  for the  weight in Eq. \eqref{sexticFreud} takes the form

\begin{align}
\label{genFrenormalizationb}
\hat{\Gamma}_m
=\langle \mathcal{S}_{m}, \mathcal{S}_{m} \rangle_{ W_{\sigma}}=\Vert \mathcal{S}_{m} \Vert^2_{ W_{\sigma}}= \sum_{k=0}^{\lfloor\frac m2\rfloor}\Psi_{k}(m)~\eta_{2m-2k}(t;\sigma),
\end{align}
where  $\Psi_{k}(m)$ is given in  Eq. \eqref{maincolsedfora}.
\noindent Eq. \eqref{genFrenormalizationb} is equivalently given by
\begin{equation*}
\hat{\Gamma}_{m}(t)=\int_{-\infty}^{\infty}\mathcal{S}^2_{m}(x,t)~
W_{\sigma}(x;t)~{\rm d}x.
\end{equation*}
\noindent
By using variable transformation  $x^2=\xi$, we have  different normalization parties as follows: 
\begin{align*}
\hat{\Gamma}_{2m}(t)&=\int_{-\infty}^{\infty}\mathcal{S}_{2m}^2(x,t) ~  W_{\sigma}(x;t){\rm d}x\\
&=2\int_{0}^{\infty}\mathcal{S}_{2m}^2(\sqrt{\xi},t)~
|\xi|^{\sigma +\frac{1}{2}}\exp\left(-[ c\xi^3+t(\xi^2 -\xi)]\right) \dfrac{1}{2\sqrt{\xi}}{\rm d}\xi\\
&=\int_{0}^{\infty}\widetilde{P}_m^2(\xi,t)s^{-\frac{1}{2}}~ |\xi|^{\sigma +\frac{1}{2}}\exp\left(-[ c\xi^3+t(\xi^2 -\xi)]\right)~{\rm d}\xi=:\widetilde{h}_m(t),
\end{align*}
and
\begin{align*}
\hat{\Gamma}_{2m+1}(t)&=\int_{-\infty}^{\infty}\mathcal{S}_{2m+1}^2(x,t) ~  W_{\sigma}(x;t)~{\rm d}x\\
&=2\int_{0}^{\infty}\mathcal{S}_{2m+1}^2(\sqrt{\xi},t)~
|\xi|^{\sigma +\frac{1}{2}}\exp\left(-[ cs^3+t(s^2 -s)]\right) \dfrac{1}{2\sqrt{\xi}}~{\rm d}\xi\\
&=\int_{0}^{\infty}\widehat{P}_n^2(\xi,t)~\xi^{\frac{1}{2}}~ |\xi|^{\sigma +\frac{1}{2}}\exp\left(-[ c\xi^3+t(\xi^2 -\xi)]\right)~{\rm d}\xi=:\widehat{h}_m(t),
\end{align*}
We now see that
\begin{align*}
\mathcal{S}_{2m}(\sqrt{\xi},t)&=(\sqrt{\xi})^{2m}+\chi(2m,t)(\sqrt{\xi})^{2m-2}+\ldots+\mathcal{S}_{2m}(0,t)\\
&=\xi^n+\widetilde{p}(m,t)\xi^{m-1}+\ldots+\widetilde{P}_m(0,t):=\widetilde{P}_m(\xi,t),
\end{align*}
and
\begin{align*}
\mathcal{S}_{2m+1}(\sqrt{\xi},t)&=(\sqrt{\xi})^{2m+1}+\chi(2m,t)(\sqrt{\xi})^{2m-1}+\ldots+k\cdot \sqrt{\xi},~~~~k\in \mathbb{R},
\\
&=\sqrt{\xi}\left(\xi^m+\widehat{p}(m,t)\xi^{m-1}+\ldots+k\right):=\sqrt{\xi} \widehat{P}_m(\xi,t).
\end{align*}
The above polynomials $\widetilde{P}_m(\xi,t)$ and $\widehat{P}_m(\xi,t)$ are recognized  as monic semiclassical Airy-type polynomials with  corresponding orthogonality weights
\begin{subequations}\label{Weights}
	\begin{align}
	\label{weight1}
	w_1(x;t)= \xi^{-\frac{1}{2}}|\xi|^{\sigma +\frac{1}{2}}\exp\left(-[ c\xi^3+t(\xi^2 -\xi)]\right),
	\end{align}
	
	\begin{align}
	\label{weight2}
	w_2(x;t)=\xi^{\frac{1}{2}}|\xi|^{\sigma +\frac{1}{2}}\exp\left(-[ c\xi^3+t(\xi^2 -\xi)]\right),
	\end{align}
\end{subequations}
both defined over $(0,\infty)$ respectively. (See \cite{refChihara} for symmetrization process and quadratic transformation).

\noindent The corresponding Hankel determinants for the weights in \eqref{Weights}  can be given by
{\small\begin{align*}
\widetilde{D}_m(t):= &{\rm det}\left(\int_{0}^\infty\xi^{i+j-\frac{1}{2}}~\xi^{\sigma +\frac{1}{2}}\exp\left(-[ c\xi^3+t(\xi^2 -\xi)]\right)~{\rm d}\xi \right)_{i,j=0}^{n-1}=\prod_{l=0}^{m-1}\widetilde{h}_l(\xi),\\
\widehat{ D}_m(t):= &{\rm det}\left(\int_{0}^\infty\xi^{i+j+\frac{1}{2}}~\xi^{\sigma +\frac{1}{2}}\exp\left(-[ c\xi^3+t(\xi^2 -\xi)]\right)~{\rm d}\xi\right)_{i,j=0}^{n-1}=\prod_{l=0}^{m-1}\widehat{h}_l(\xi)
\end{align*}}
respectively. Hence,
\begin{equation*}
\Delta_n(t)=\prod_{j=0}^{n-1}\Gamma_j(t)=\begin{cases}
\widetilde{D}_{k+1}\widehat{ D}_k&~~ \text{$n=2k+1$,}\\
\widetilde{D}_k\widehat{ D}_k& ~~\text{$n=2k$.}
\end{cases}
\end{equation*}
\noindent It is good to mention here that investigation of asymptotics of the  Hankel determinants when  $n$ is large    has been an interesting subject for many years; for instance, for Gaussian weight is studied in Chen et al. in \cite{LC2}. See also the monograph by Szeg$\ddot{o}$ \cite{refSzego} as we will not address this as it goes beyond the scope of  the paper.
\subsection{Nonlinear recursion relation}
In this section, we explore certain nonlinear recurrence relations associated with the semi-classical weight given in \eqref{sexticFreud}.
\begin{theorem}{\label{Wronskianch3}}{ 
		For the semiclassical weight in \eqref{orthosextic}, the recurrence coefficient $\gamma_n(t;\sigma)$ 	fullfill  the following difference relations
		\begin{align}\label{nonlinerreucrequnforFreudAnewAx}
		6c\left[
		\gamma_{n}\left( \Xi_{n-1} +\Xi_{n} +\Xi_{n+1}\right) + \gamma_{n-1} \gamma_{n}  \gamma_{n+1} \right] + 4t\Xi_{n}-2t\gamma_{n}
		=	n+ (2\sigma +1) \Omega_{n},
		\end{align}
		with initial conditions given by
		\begin{align}\label{intcond}
		\begin{cases}
		&
		\gamma_1(t;\sigma)
		= \dfrac{ \Vert x^2\Vert_{t}^{2}}{\Vert 1\Vert_{t}^{2}} =\dfrac{\eta_2(t;\sigma)}{\eta_0(t;\sigma)}
		=\frac{\int_{-\infty}^{\infty} x^2  W_{\sigma}(x;t)\mathrm{d}x}{\int_{-\infty}^{\infty}  W_{\sigma}(x;t)~\mathrm{d} x},
		\\&
		\gamma_{0}=0,
		\end{cases}
		\end{align} 
		where $\Xi_n$ and $\Omega_{n}$ are, respectively, given by
		\begin{align}\label{greatqn}
		\Xi_n= \gamma_n(t;\sigma) \left[\gamma_{n-1}(t;\sigma)+\gamma_n(t;\sigma) +\gamma_{n+1}(t;\sigma)\right],
		\end{align}
		and
		\begin{align}\label{OmegaN}
		\Omega_n= \dfrac{1-(-1)^n}{2} = \begin{cases}
		1,~~ \text{for}\quad \text{$n$ is odd}
		\\ 0,~~ \text{for}\quad \text{$n$ is even}.
		\end{cases}
		\end{align}
}\end{theorem}
\begin{proof}\begin{itemize}
		\item[(i)] 
		Applying similar procedure due to Freud as given in \cite[Section 2]{van2007discrete} (see also \cite{refnevai1986geza}),  let's consider the following integral
		\begin{equation}\label{mainIn}
			\mathbb{J}_n = \frac{1}{\hat{\Gamma}_{n}} \displaystyle \int_{-\infty}^{\infty} \left[ \mathcal{S}_{n}(x;t)~\mathcal{S}_{n-1}(x;t)\right]^{'}~W_{\sigma}(x;t) ~\mathrm{d} x ,
			\end{equation}
		where 	$\hat{\Gamma}_{n}$ is  given in  \eqref{genFrenormalizationb}. 
		Eq. \eqref{mainIn} is equivalently given by
		{\small	\begin{align}\label{Freudintegralz}
			\mathbb{J}_n &
			=  	\frac{1}{\hat{\Gamma}_{n}} 
			\Big[\langle  S^{'}_n,\mathcal{S}_{n-1}\rangle_{W_{\sigma}}  +
			\langle \mathcal{S}_{n},\mathsf{S}^{'}_{n-1}\rangle_{W_{\sigma}}\Big]
			\nonumber \\&	=  \frac{1}{\hat{\Gamma}_{n}} \displaystyle \int_{-\infty}^{\infty}  \left(nx^{n-1}+V_{n-2}\right)~\mathcal{S}_{n-1}(x;t)~W_{\sigma}(x;t) ~\mathrm{d} x 
			= \frac{\hat{\Gamma}_{n-1}}{\hat{\Gamma}_{n}} n,
			\end{align}}
		where $V_{n-2} \in \mathbb{P}_{n-2}.$ 
		We also see that by evaluating Eq. \eqref{mainIn} using technique of integration, we arrive at		
		{\small\begin{align}\label{coolmainAA}
			\nonumber\mathbb{I}_n\hat{\Gamma}_{n}&= 
			\left[ \mathcal{S}_{n}(x;t)~\mathcal{S}_{n-1}(x;t) W_{\sigma}(x;t)\right]_{-\infty}^{\infty}
			-   \displaystyle \int_{-\infty}^{\infty}  \mathcal{S}_{n}(x;t)~\mathcal{S}_{n-1}(x;t)~W^{'}_{\sigma}(x;t)\medspace \mathrm{d} x 
			\nonumber\\& = 
			-  (2\sigma +1) \displaystyle \int_{-\infty}^{\infty}  \dfrac{\mathcal{S}_{n}(x;t)~\mathcal{S}_{n-1}(x;t)}{x}~W_{\sigma}(x;t)\medspace \mathrm{d} x 
			+ 6c \displaystyle \int_{-\infty}^{\infty}  x^5\mathcal{S}_{n}(x;t)~\mathcal{S}_{n-1}(x;t)~W_{\sigma}(x;t)\medspace \mathrm{d} x 
			\nonumber
			\\& \qquad 	+4t \displaystyle \int_{-\infty}^{\infty} x^3 \mathcal{S}_{n}(x;t)~\mathcal{S}_{n-1}(x;t)~ ~W_{\sigma}(x;t)\medspace \mathrm{d} x 
			- 2t \displaystyle \int_{-\infty}^{\infty}  x\mathcal{S}_{n}(x;t)~\mathcal{S}_{n-1}(x;t)~ ~W_{\sigma}(x;t)\medspace \mathrm{d} x ,		\end{align}}
		in consideration of the fact that {\small$\Big[\mathcal{S}_{n}(x;t)~\mathcal{S}_{n-1}(x;t)  W_{\sigma}(x;t)\Big]_{-\infty}^{\infty}=0$} as the weight \eqref{orthosextic} vanishes at the  boundary terms  when $x\rightarrow \pm \infty$ due to symmetry property of the weight 	$W_{\sigma}$; hence it follows that	
		{\small	\begin{subequations}\label{patAS}
				{\small\begin{align}\label{patAS1}
					\displaystyle \int_{-\infty}^{\infty}  \mathcal{S}_{n}(x;t)~ \dfrac{1}{x} ~\mathcal{S}_{n-1}(x;t)~ W_{\sigma}(x;t)\medspace \mathrm{d} x  =0,
					\end{align}}
				for $n$ is even and, when $n$ is odd, we have that
				\begin{align}\label{patAS2}
				\displaystyle \int_{-\infty}^{\infty}\mathcal{S}_{n-1}(x;t) ~  \dfrac{\mathcal{S}_{n}(x;t)}{x} ~ ~W_{\sigma}(x;t)\medspace \mathrm{d} x  =\hat{\Gamma}_{n-1},
				\end{align}
				as
				$\dfrac{\mathcal{S}_{n}(x;t)}{x}$ is a polynomial of degree $n-1$.
				Thus, we have
				\begin{align}\label{patAS2ZZ}
				\displaystyle \int_{-\infty}^{\infty} ~  \dfrac{\mathcal{S}_{n-1}(x;t) \mathcal{S}_{n}(x;t)}{x}  ~W_{\sigma}(x;t)\medspace \mathrm{d} x  =\Omega_n \hat{\Gamma}_{n-1},
				\end{align}
		\end{subequations}}
where  {\small$\Omega_n$} is given in \eqref{OmegaN}. Let's us employ the following iterated recurrence relation from Eq. \eqref{symmetricrecur} to obtain
		\begin{subequations}\label{recurrencidnetities}
			\begin{align}
			x^5 \mathcal{S}_n(x;t)& = \mathcal{S}_{n+5}(x;t) + 
			\left( \gamma_{n} + \gamma_{n+1}  + \gamma_{n+2} + \gamma_{n+3} +\gamma_{n+4}\right)\mathcal{S}_{n+3}(x;t) 
			\nonumber \\&\quad
			+ \left[ \gamma_{n}\left( \Xi_{n-1} +\Xi_{n}+ \Xi_{n+1} \right) + \gamma_{n-1}\gamma_{n}  \gamma_{n+1}  \right]\mathcal{S}_{n+1}(x;t) 
			\nonumber \\&\quad				
			+ 			\left[\gamma_{n}\gamma_{n-2}  \Xi_{n-1} + \gamma_{n-2} \gamma_{n-1}\gamma_{n}\gamma_{n+1} + \gamma_n \gamma_{n-1} \gamma_{n-2}\gamma_{n-3} \right]\mathcal{S}_{n-3}(x;t)
			\nonumber \\&\quad
			+\left(\gamma_{n}\gamma_{n-1}\gamma_{n-2}\gamma_{n-3}\gamma_{n-4}\right)
			\mathcal{S}_{n-5}(x;t), \label{recurrence5}
			\end{align}
			
			\begin{align}
			x^4 \mathcal{S}_n(x;t)& = \mathcal{S}_{n+4}(x;t) + (  \gamma_{n}+ \gamma_{n+1} + \gamma_{n+2} +\gamma_{n+3})  \mathcal{S}_{n+2}(x;t)
			\nonumber\\& \quad 
			+ \big[\gamma_{n} ( \gamma_{n-1} + \gamma_{n}  + \gamma_{n+1})    +\gamma_{n+1} ( \gamma_{n} + \gamma_{n+1}  + \gamma_{n+2} )\big]  \mathcal{S}_{n}(x;t)
			\nonumber\\& \quad  +\gamma_{n}  \gamma_{n-1}(  \gamma_{n-2} +\gamma_{n-1} + \gamma_{n}  + \gamma_{n+1} ) \mathcal{S}_{n-2}(x;t)
			+ ( \gamma_{n}  \gamma_{n-1} \gamma_{n-2} \gamma_{n-3} )  \mathcal{S}_{n-4}(x;t),\label{recurrence3}
			\end{align}
			\begin{align}
			x^3\mathcal{S}_{n}(x;t)& = \left(\gamma_{n} +\gamma_{n+1}+\gamma_{n+2}\right)	\mathcal{S}_{n+1}(x;t) +\mathcal{S}_{n+3}(x;t) 
			\nonumber \\&\quad
			+ \gamma_{n} \gamma_{n-1}~\gamma_{n-2}\mathcal{S}_{n-3}(x;t)	+ \gamma_{n} \left(\gamma_{n-1} +\gamma_{n}+\gamma_{n+1}\right)
			\mathcal{S}_{n-1}(x;t), \label{recurrence2}
			\end{align}
			
			\begin{align}
			x^2\mathcal{S}_{n}(x;t) &= \left( \gamma_{n} + \gamma_{n+1}  \right)\mathcal{S}_{n}(x;t) + \gamma_{n} \gamma_{n-1}\mathcal{S}_{n-2}(x;t)+\mathcal{S}_{n+2}(x;t).
			\end{align}
		\end{subequations}
\noindent By using the identities \eqref{recurrencidnetities} and Eq. \eqref{eq:Pearson} for the weight \eqref{sexticFreud} together with Eqs. \eqref{patAS} into \eqref{coolmainAA}, we obtain 				
		\begin{align}\label{nonlinerreucrequnforFreud}
		\nonumber
		n\hat{\Gamma}_{n-1}=\mathbb{I}_n\hat{\Gamma}_{n}&= 6c\left[\left(\gamma_{n} +\gamma_{n-1}\right)\Xi_n + 
		\left( \gamma_{n} \Xi_{n+1}+ \gamma_{n}  \gamma_{n-1} \gamma_{n-2} \right)\right]\hat{\Gamma}_{n-1}
		\\& \quad -2t\gamma_n \hat{\Gamma}_{n-1} -  (2\sigma +1) \Omega_{n} \hat{\Gamma}_{n-1}+4t\left[\gamma_{n} \left(\gamma_{n-1} +\gamma_{n}+\gamma_{n+1}\right)\right]\hat{\Gamma}_{n-1},
		\end{align}
		which simplifies, using the fact that {\small$\hat{\Gamma}_{n-1}\neq 0$}, to
		\begin{align}\label{nonlinerreucrequnforFreudA}
		\nonumber
		n+ (2\sigma +1) \Omega_{n} &= 6c\left[\left(\gamma_{n} +\gamma_{n-1}\right)\Xi_n + 
		\left( \gamma_{n} \Xi_{n+1}+ \gamma_{n}  \gamma_{n-1} \gamma_{n-2} \right)\right]
		\\&
		+ 4t\left[\gamma_{n} \left(\gamma_{n-1} +\gamma_{n}+\gamma_{n+1}\right)\right] -2t\gamma_n,
		\end{align}
		where  {\small$\Omega_{n}$} is given in \eqref{OmegaN}.  Note that Eq. \eqref{nonlinerreucrequnforFreud} and Eq. \eqref{Freudintegralz} yield Eq. \eqref{nonlinerreucrequnforFreudAnewAx}.
	\end{itemize}
\end{proof}
\begin{remark}
	Quite similar non-linear  discrete equations  like Eq. \eqref{nonlinerreucrequnforFreudA} can be obtained in  \cite[Eq. (23), p. 5]{freud1976coefficients} and we also refer to \cite{aptekarev1997toda,clarkson2016generalized,van2017orthogonal}.
\end{remark}
\noindent
The following result gives the differential-recurrence relation for the weight \eqref{sexticFreud}.
\begin{theorem}\label{TodaLike}
For the semiclassical weight in \eqref{orthosextic},  the  coefficients $\gamma_n(t;\sigma)$ obey Toda-type formulation
	\begin{align}\label{diffrecursextic}
	\dfrac{\,d \gamma_n}{\,d{t}}
	= \gamma_n \big[\left(\gamma_{n+1}-\Xi_{n+1} \right) - \left(\gamma_{n-1} -\Xi_{n-1}\right)\big],
	\end{align}
	where $\Xi_n$ is given in  Eq. \eqref{greatqn}.
\end{theorem}
\begin{proof}
	In order to prove this result, we first differentiate the normalization constant {\small$\hat{\Gamma}_{n}(t)$} with respect to $t$ as	
	\begin{align}\label{DiffereintalnormalisedSextt}
	\dfrac{\,d \hat{\Gamma}_{n}}{\,d{t}}
	& = 2\langle\dfrac{\,d \mathcal{S}_{n}}{\,d{t}}, \mathcal{S}_{n} \rangle_{ W_{\sigma}}  +  \langle (x^2-x^4)\mathsf{S}_{n},\mathsf{S}_{n} \rangle_{ W_{\sigma}},
	\nonumber \\ &
	= 2\int_{-\infty}^{\infty} \dfrac{\,d \mathcal{S}_{n}(x;t)}{\,d{t}}~\mathcal{S}_{n}(x;t)~W_{\sigma}(x;t)~\mathrm{d} x  +	\int_{-\infty}^{\infty} x^2~\mathsf{S}^2_{n}(x;t)~W_{\sigma}(x;t)~\mathrm{d} x 
	\\&  - 	\int_{-\infty}^{\infty} x^4~\mathsf{S}^2_{n}(x;t)~W_{\sigma}(x;t)~\mathrm{d} x .
	\nonumber	\end{align}
	We see from Eq. \eqref{DiffereintalnormalisedSextt} that  the first integral vanishes by orthogonality  as
	$\dfrac{\,d\mathcal{S}_{n}}{\,d{t}}\in \mathcal{P}_{n-1}$.
Using the recursive  relation in Eq. \eqref{symmetricrecur} and orthogonality fact, we  now have
	\begin{align}\label{DiffereintalnormalisedSexttb}
	\dfrac{\,d}{\,d{t}}\hat{\Gamma}_{n}
	=		\left(\gamma_n + \gamma_{n+1}\right)\hat{\Gamma}_{n}  - \left(\Xi_n +\Xi_{n+1} \right)\hat{\Gamma}_{n}
	=\left[\left(\gamma_n -\Xi_n\right)  + \left(\gamma_{n+1}-\Xi_{n+1} \right)\right]\hat{\Gamma}_{n},
	\end{align}
	Besides, if we  differentiate Eq. \eqref{betafromularecur} with respect to $t$, we obtain
	\begin{align}\label{coolA}
	\dfrac{\,d}{\,d{t}}\gamma_{n}
	=
	\dfrac{\,d}{\,d{t}}\left(\dfrac{\hat{\Gamma}_{n}}{\hat{\Gamma}_{n-1}}\right)
	=\gamma_n \left[
	\dfrac{\,d}{\,d{t}}\ln \hat{\Gamma}_{n} -
	\dfrac{\,d}{\,d{t}}
	\ln \hat{\Gamma}_{n-1}\right]
	=\gamma_n \Big[
	\left(\gamma_{n+1}-\gamma_{n-1} \right) - \left[\Xi_{n+1} -\Xi_{n-1}\right]\Big],
	\end{align}
	and substituting  Eq. \eqref{DiffereintalnormalisedSexttb} into  \eqref{coolA} leads to the required result.
\end{proof}
The following result presents  nonlinear differential-recurrence relation of high order associated with the weight \eqref{orthosextic}; we quote ideas of the proof from \cite{refkelilAppadu2020}.
\begin{theorem}
	The  coefficients $\gamma_n(t;\sigma)$ for the weight in  Eq.  \eqref{sexticFreud} fulfills the following nonlinear differential-recurrence equation
	{\small
		\begin{align}\label{secondorderrecur}
		\begin{cases}
		\dfrac{\,d^2 \gamma_n}{\,d{t}}
		&=  \frac{1}{6c} \left[n+(2\sigma +1) \Omega_n - \vartheta(t)\right] + \left( -\gamma_{n-1}-\gamma_{n+1} \right) \gamma_{n}^{4}
		\nonumber\\&\qquad 
		+ \left( -\gamma_{n-2}\gamma_{n-1}-\gamma^{2}_{n-1}-6\gamma_{n-1}\gamma_{n+1}-\gamma^{2}_{n+1}-\gamma_{n+1}\gamma_{n+2}
		+2\gamma_{n-1}+2\gamma_{n+1} \right) \gamma^{3}_{n}
		\nonumber\\&\qquad 
		+
		\Bigg( \gamma_{n-3}\gamma_{n-2}\gamma_{n-1}+\gamma^{2}_{n-2}\gamma_{n-1}+2\gamma_{n-2}\gamma^{2}_{n-1}-4\gamma_{n-2}\gamma_{n-1}\gamma_{n+1}
		+\gamma^{3}_{n-1}
		-5\gamma^{2}_{n-1}\gamma_{n+1}	-4\gamma_{n-1}\gamma_{n+1}\gamma_{n+2}
		\nonumber\\& \qquad \qquad
		-5\gamma_{n-1}\gamma^{2}_{n+1} 	+\gamma^{3}_{n+1} 	+2\gamma^{2}_{n+1}\gamma_{n+2}  +\gamma_{n+1}\gamma^{2}_{n+2}+\gamma_{n+1}\gamma_{n+2}\gamma_{n+3}+8\gamma_{n-1}\gamma_{n+1}-\gamma_{n-1}-\gamma_{n+1} 
		\Bigg) \gamma^{2}_{n}
		\nonumber\\&\qquad +
		\Bigg(
		\gamma_{n-4}\gamma_{n-3}\gamma_{n-2}\gamma_{n-1}+\gamma^{2}_{n-3}\gamma_{n-2}\gamma_{n-1}
		+2\gamma_{n-3}\gamma^{2}_{n-2}\gamma_{n-1}+2\gamma_{n-3}\gamma_{n-2}\gamma^{2}_{n-1}+\gamma^{3}_{n-2}\gamma_{n-1}
		\nonumber\\& \qquad\qquad
		+3\gamma^{2}_{n-2}\gamma^{2}_{n-1}+3\gamma_{n-2}\gamma^{3}_{n-1}-2\gamma_{n-2}\gamma_{n-1}\gamma^{2}_{n+1}
		-2\gamma_{n-2}\gamma_{n-1}\gamma_{n+1}\gamma_{n+2}+\gamma^{4}_{n-1}-2\gamma^{2}_{n-1}\gamma^{2}_{n+1}
		\nonumber\\& \qquad \qquad 
		-2\gamma^{2}_{n-1}\gamma_{n+1}\gamma_{n+2}+\gamma^{4}_{n+1}+3\gamma^{3}_{n+1}\gamma_{n+2}+3\gamma^{2}_{n+1}\gamma^{2}_{n+2}
		+2\gamma^{2}_{n+1}\gamma_{n+2}\gamma_{n+3}+\gamma_{n+1}\gamma^{3}_{n+2}	
		\nonumber\\& \qquad \qquad 
		+2\gamma_{n+1}\gamma^{2}_{n+2}\gamma_{n+3}+\gamma_{n+1}\gamma_{n+2}\gamma^{2}_{n+3}
		-2\gamma^{2}_{n-2}\gamma_{n-1}
		+\gamma_{n+1}\gamma_{n+2}\gamma_{n+3}\gamma_{n+4}-2\gamma_{n-3}\gamma_{n-2}\gamma_{n-1}
		\nonumber\\& \qquad \qquad
		-4\gamma_{n-2}\gamma^{2}_{n-1}
		+2\gamma_{n-2}\gamma_{n-1}\gamma_{n+1} 	-2\gamma^{3}_{n-1}
		+2\gamma^{2}_{n-1}\gamma_{n+1}+2\gamma_{n-1}\gamma^{2}_{n+1}
		+2\gamma_{n-1}\gamma_{n+1}\gamma_{n+2}-2\gamma^{3}_{n+1}
		\nonumber \\& \qquad \qquad 
		-4\gamma^{2}_{n+1}\gamma_{n+2} 	-\gamma^2_n	-2\gamma_{n+1}\gamma^{2}_{n+2}-2\gamma_{n+1}\gamma_{n+2}\gamma_{n+3}-2\gamma_{n-1}\gamma_{n+1} 
		-2\gamma_n\gamma_{n-1} - 2\gamma_n \gamma_{n+1}-\gamma_{n+1}\gamma_{n-1}
		\Bigg)\gamma_{n},
		\\& 
		\vartheta(t) = 4t\Xi_n-2\gamma_n t = 2t\gamma_n \left[2(\gamma_{n-1}+\gamma_{n}+\gamma_{n+1})-1\right].
		\end{cases}
		\end{align}
	}
	where {\small$\Omega_{n}$} and {\small$\Xi_n$} are given in \eqref{OmegaN} and \eqref{greatqn} respectively.
\end{theorem}
\begin{proof}
	For the proof, we refer similar ideas in \cite{refkelilAppadu2020}. 
\end{proof}
\subsection{Differential-Recurrence relation}\label{differdiff}
\noindent Chen and Feigin \cite{chen2006painleve} obtained ladder operators for  a semiclassical weight $\widetilde{w}(x)|x-t|^{\varTheta},$ where $ x, \varTheta, t\in\mathbb{R}$ and $\widetilde{w}(x)$ is classical weight function. 
	In Filipuk et al. \cite{filipuk2012discrete}, it is shown that the recurrence coefficients for the quartic Freud weight {\small$|x|^{2\alpha+1}{\rm e}^{-x^{4}+tx^{2}}, x,t\in\mathbb{R}, \medspace \alpha>-1$} are related to the solutions of the Painlev\'{e} IV and the first discrete Painlev\'{e} equation. Clarkson et al. \cite{clarkson2016generalized} provided a systematic study on Freud weights and some generalized work for  \cite{chen2006painleve}.
\begin{lemma}{\cite{refkelilAppadu2020}}
	\label{lem5.1}
	The monic orthogonal polynomials $P_{n}(x;t)$ with respect to the semiclassical Freud-type weight \eqref{sexticFreud}
	\begin{equation*}
		w_{\alpha}(x)=|x|^{\alpha}w_{0}(x),
		\end{equation*}
	where
	\begin{equation*}
		w_{0}(x):={\rm e}^{-v_{0}(x) }~~~ {\rm with} ~~~v_{0}(x):=cx^6+t(x^4 -x^2).
		\end{equation*}
	on $\mathbb{R}$ satisfy the differential-difference recurrence relation
	{\small\begin{subequations}
			\begin{equation*}
			P'_{n}(x)=\gamma_{n}(t)\mathcal{A}_{n}(x)P_{n-1}(x)-\mathcal{B}_{n}(x)P_{n}(x),
			\end{equation*}
			where
			\begin{equation}\label{CoeffAna}
			\mathcal{A}_{n}(x):=\frac{1}{\Gamma_{n}}\int_{-\infty}^{\infty}\frac{v_{0}'(x)-v_{0}'(\tau)}{x-\tau}~P^{2}_{n}(\tau)~w(\tau)\,d\tau,\\
			\end{equation}
			\begin{equation}\label{Coeffbn}
			\mathcal{B}_{n}(x):=\frac{1}{\Gamma_{n-1}}\int_{-\infty}^{\infty}\frac{v_{0}'(x)-v_{0}'(\tau)}{x-\tau}~P_{n}(\tau)~P_{n-1}(\tau)~w(\tau)\,d\tau+\frac{\alpha\left[1-(-1)^{n}\right]}{2x}.
			\end{equation}
	\end{subequations}}
\end{lemma}
\begin{proof}
	For the proof, we refer to \cite{refkelilAppadu2020}. See also similar works in \cite{chen1997ladder}. 
\end{proof}
\begin{lemma}\label{lem5.2}
		$\mathcal{A}_{n}(z)$ and $\mathcal{B}_{n}(z)$ defined by Lemma \ref{lem5.1} satisfy the relation:
		\begin{equation}
		\label{AnBn}
		\mathcal{A}_{n}(z)=\frac{v_{0}'(z)}{z}+\frac{\mathcal{B}_{n}(z)+\mathcal{B}_{n+1}(z)}{z}-\frac{\alpha}{z^{2}}.
		\end{equation}
\end{lemma}
\begin{proof}
	Be the definition of {\small$\mathcal{A}_{n}(z)$}, we rewrite it as
	{\small\begin{equation*}
		\begin{split}
		\mathcal{A}_{n}(z)=&\frac{1}{z\Gamma_{n}}\left\{\int_{-\infty}^{\infty}\frac{v_{0}'(z)-v_{0}'(\tau)}{z-\tau}yP_{n}^{2}(\tau)w(\tau)\mathrm{d}\tau+\int_{-\infty}^{\infty}\left[v_{0}'(z)-v_{0}'(\tau)\right]P_{n}^{2}(\tau)w(\tau)\mathrm{d}\tau\right\}\\
		=&\frac{1}{z\Gamma_{n}}\bigg\{\int_{-\infty}^{\infty}\frac{v_{0}'(z)-v_{0}'(\tau)}{z-\tau}\left[P_{n+1}(\tau)+\gamma_{n}P_{n-1}(\tau)\right]P_{n}(\tau)w(\tau)\mathrm{d}\tau+v_{0}'(z)\Gamma_{n}
		\\&\qquad \qquad 
		-\int_{-\infty}^{\infty}P_{n}^{2}(\tau)\left[\frac{\alpha}{\tau}w(\tau)-w'(\tau)\right]\mathrm{d}\tau\bigg\}\\
		=&\frac{1}{z}\left\{\mathcal{B}_{n+1}(z)-\frac{\alpha}{2z}\left[1-(-1)^{n+1}\right]+\mathcal{B}_{n}(z)-\frac{\alpha}{2z}\left[1-(-1)^{n}\right]\right\}+\frac{v_{0}'(z)}{z},\\
		=&\frac{\mathcal{B}_{n}(z)+\mathcal{B}_{n+1}(z)}{z}-\frac{\alpha}{z^{2}}+\frac{v_{0}'(z)}{z},
		\end{split}
		\end{equation*}}
	which completes the proof.
\end{proof}
\begin{lemma}{\cite[Chapter 3]{refIsmaila}.}
	The functions $\mathcal{A}_n(z)$, $\mathcal{B}_n(z)$ and $\sum_{k=0}^{n-1}\mathcal{A}_k(z)$ satisfy the identity
	\begin{align}\label{compCC}
	\mathcal{B}^2_{n}(z) +v{'}(z)\mathcal{B}_{n}(z) +
	\sum_{k=0}^{n-1}\mathcal{A}_k(z)=\gamma_n\mathcal{A}_n(z)\mathcal{A}_{n-1}(z).
	\end{align}
\end{lemma}
We, next, apply the ladder coefficients to the case of perturbed Freud weight as follows.
\subsubsection{Ladder operator relations for the  weight \eqref{sexticFreud}}\label{sc}
For the perturbed Freud-type weight \eqref{sexticFreud}, 
\begin{align}\label{nusextiFreud1a}
	v(x) = -\ln  W_{\sigma}(x;t)  = -(2\sigma +1)\ln |x| +cx^6+t(x^4-x^2),~~x\in \mathbb{R},
	\end{align}
we have 
\begin{align*}
	v{'}(x) = -\frac{(2\sigma +1)}{x} +6cx^5+t(4x^3-2x),
	\end{align*}
and hence
\begin{align*}
	\dfrac{v{'}(x) -v{'}(\tau)}{x-\tau}=\dfrac{2\sigma +1}{x\tau}+ 6c \{x^4+x^3\tau+x^2\tau^2+x\tau^3+\tau^4\}+4t(x^2+x\tau+\tau^2)-2t.
	\end{align*}
\begin{theorem}\label{TheoremDEforSexticFreud}
	The monic orthogonal polynomials $\mathcal{S}_{n}(x;t)$ with respect to  the weight in \eqref{sexticFreud}
	defined  on {\small$\mathbb{R}$} obey the relation
	\begin{align*}
		\mathsf{S}'_n(x;t) = \gamma_n(t)\mathcal{A}_n(x;t)\mathcal{S}_{n-1}(x;t) -\mathcal{B}_n(x;t)\mathcal{S}_{n}(x;t)
		\end{align*}
	where 
	{\small\begin{subequations}\label{CoeffpertrubeFreudone}
			\begin{align}
			\mathcal{A}_n(x;t) = 6cx^4 +6c(\gamma_n+\gamma_{n+1})x^2+ 6c\left(\Xi_{n+1} +\Xi_{n}\right)  +4tx^2 +4t(\gamma_n+\gamma_{n+1})  -2t,
			\end{align}
			\begin{align}
			\mathcal{B}_n(x;t) = \left(\dfrac{2\sigma +1}{x}\right) \Omega_n+ 6c\gamma_{n}x^3 +6c \Xi_{n}x +4tx\gamma_{n},
			\end{align}
	\end{subequations}}
	where the expressions  {\small$\Xi_n$} and {\small$\Omega_n$} are  given   in \eqref{greatqn}  and   \eqref{OmegaN} respectively.
\end{theorem}
\begin{proof}
	From \eqref{CoeffAna}, we obtain
	\begin{align}\label{CoeffAnGenSextic}
	\mathcal{A}_n(x;t)&= \dfrac{1}{\hat{\Gamma}_{n}}\int_{\mathbb{R}}\mathsf{S}^2_{n}(\tau) \left(\dfrac{v{'}(x) -v{'}(\tau)}{x-\tau}\right)W_{\sigma}(\tau;t)\mathrm{d}\tau
	\nonumber\\&= \dfrac{1}{\hat{\Gamma}_{n}}\int_{\mathbb{R}}\mathsf{S}^2_{n}(\tau)\Big(\dfrac{2\sigma +1}{x\tau}+ 6c \{x^4+x^3\tau+x^2\tau^2+x\tau^3+\tau^4\}	+4t(x^2+x\tau+\tau^2)-2t\Big) W_{\sigma}(\tau;t)\mathrm{d}\tau
	\nonumber\\&=
	6cx^4 +6c(\gamma_n+\gamma_{n+1})x^2+ 6c\left(\Xi_{n+1} +\Xi_{n}\right)  +4tx^2 +4t(\gamma_n+\gamma_{n+1})  -2t,  
	\end{align}
	and the integral  in \eqref{CoeffAnGenSextic}  vanishes due to symmetry of  $W_{\sigma}$.
	
	\noindent Besides,  by using Eq. \eqref{Coeffbn}, orthogonality and Eq. \eqref{symmetricrecur}, we  have that
	\begin{align}\label{CoeffBnGenSextic}
	\nonumber
	\mathcal{B}_n(x;t)
	&= \dfrac{1}{\hat{\Gamma}_{n-1}}\int_{\mathbb{R}}\mathcal{S}_{n}(\tau) \mathcal{S}_{n-1}(\tau) \Big(\dfrac{2\sigma +1}{xy}+ 6c \{x^4+x^3y+x^2y^2+xy^3+y^4\}
	\\&\qquad \qquad\qquad\qquad\qquad\qquad
	+4t(x^2+xy+y^2)-2t\Big) W_{\sigma}(y;t)\,d{y}
	\nonumber\\&=   6c\gamma_{n}x^3 +6c \Xi_{n}x +4tx\gamma_{n} + \left(\dfrac{2\sigma +1}{x}\right) \Omega_n,
	\end{align}
	where  $\Xi_n$ and $\Omega_n$ are given respectively in \eqref{greatqn} and \eqref{OmegaN}.
\end{proof}
\begin{remark}
	It is good to mention that there is a similar result in \cite{CJ2020} for differential-recurrence relation for sextic Freud-type weight; whereas our considered weight in Eq. \eqref{sexticFreud} can be perceived as generalized  measure deformation  using $d\mu(x;t)=e^{t(x^4-x^2)}d\mu(x;0),$  For a similar procedure, one can see \cite{refHanCHen} where the  authors used classical measure deformation	via $d\mu(x;t)=e^{tx^2}d\mu(x;0)$ for Laguerre-type weight. 
\end{remark}

\subsection{Shohat's quasi-orthogonality method}\label{Shohat's method}
\noindent
Shohat \cite{shohat1939differential} studied a strategy  using quasi-orthogonality, to find differential-difference relation for
a general semiclassical weight function.    Bonan, Freud, Mhaskar and Nevai are renowned experts who used this method in their work \cite{refnevai1986geza}.
The idea of quasi-orthogonality is well articulated in \cite{refMaroni1987,refDriverJordaan,shohat1939differential}).  Our goal in  this section is to apply this method to the case of perturbed Freud-type weight in \eqref{sexticFreud} \cite[Section 4.5]{clarkson2016generalized}.  
Following the ideas in  \cite{refnevai1986geza}, we notice that monic 
perturbed Freud-type polynomials obey  quasi-orthogonality of order {\small$m=7$} and  therefore
\begin{align}\label{quasi}
x\dfrac{\,d \mathcal{S}_n(x;\tau)}{\,\mathrm{d}x}=  \sum_{k=n-6}^{n}\mathfrak{u}_{n,k} ~\mathcal{S}_k(x;\tau),
\end{align}
where the expression  $\mathfrak{u}_{n,k}$ is obtained by
\begin{align}\label{Ccoeff}
\mathfrak{u}_{n,k} &= \frac{1}{\Gamma_{k}}  \int_{-\infty}^{\infty}   x~
\dfrac{\,d \mathcal{S}_n}{\mathrm{d}x}
(x;\tau) ~ \mathcal{S}_k(x;t)~W_{\sigma}(x;\tau)\,\mathrm{d}x,
\end{align}
with {\small$n-6\leq k\leq n$} and {\small$\Gamma_k\neq 0$}.
By employing integration techniques, for $n-6 \leq j\leq n-1$, we have
{\small\begin{align}\nonumber
	\Gamma_k ~\mathfrak{u}_{n,k} &=\Big[ x~\mathcal{S}_k(x;t)~ \mathcal{S}_n(x;t)~ W_{\sigma}(x;t)\Big]_{-\infty}^{\infty}  
	- \int_{-\infty}^\infty \dfrac{\,d}{\mathrm{d} x }\left(x \mathcal{S}_k (x;t)W_{\sigma}(x;t)\right)~\mathcal{S}_n(x;t)  \,\mathrm{d}x 
	\nonumber\\
	&=- \int_{-\infty}^{\infty}   \left[\mathcal{S}_n(x;t)~ \mathcal{S}_k(x;t)  
	+ x~\mathcal{S}_n(x;t) ~\dfrac{\mathcal{S}_k}{x}(x;t)\right] ~W_{\sigma}(x;t) \,\mathrm{d}x
	\nonumber\\&	- \int_{-\infty}^{\infty}   x\mathcal{S}_n(x;t)~ \mathcal{S}_j(x;t)~ \dfrac{\,d{W_{\sigma}(x,t)}}{\mathrm{d} x }(x;t) \,\mathrm{d}x,
	\nonumber\\
	\\
	&=    -\int_{-\infty}^{\infty} {\mathcal{S}_n(x;t) ~\mathcal{S}_j(x;t)~ \left( -6cx^{6}-4tx^4+2tx^2+ 2\sigma +1\right) }   W_{\sigma}(x;t)\,\mathrm{d}x\nonumber
	\\&= 	\int_{-\infty}^{\infty} \big(6cx^{6}+4tx^4-2tx^2- (2\sigma +1)\big) ~\mathcal{S}_n(x;t) ~\mathcal{S}_j(x;t) ~W_{\sigma}(x;t) \,\mathrm{d}x,\label{maineq}
	\end{align}}
since\[ x~
	\dfrac{\,d{W_{\sigma}(x,t)}}{\mathrm{d} x }	=\big[-6cx^{6}-4tx^4+2tx^2+ 2\sigma +1\big]W_{\sigma}(x;t).\]

\noindent
The following relations  follow from iterating the recurrence given in Eq. \eqref{symmetricrecur}:
\begin{subequations}\label{recurrence13a}
	{\small\begin{align}
		\nonumber
		x^6\mathcal{S}_{n}&(x;t)
		= \mathcal{S}_{n+6}(x;t) +  \left( \gamma_{n} + \gamma_{n+1}  + \gamma_{n+2} +\gamma_{n+3} + \gamma_{n+4}+\gamma_{n+5}\right)\mathcal{S}_{n+4}(x;t) 
		\nonumber \\&\medspace \nonumber
		+ \bigg[ \gamma_{n+3}\left( \gamma_{n} + \gamma_{n+1}  + \gamma_{n+2} + \gamma_{n+3}+ \gamma_{n+4}\right)
		+\gamma_{n+2}\left( \gamma_{n} + \gamma_{n+1}  + \gamma_{n+2} + \gamma_{n+3}\right)  +\Xi_n +\Xi_{n+1}
		\bigg]\mathcal{S}_{n+2}(x;t)
		\nonumber \\&\medspace \nonumber
		+ \bigg[\gamma_{n+1}\gamma_{n+2}\left[\gamma_{n} + \gamma_{n+1}  + \gamma_{n+2} + \gamma_{n+3}\right]+\left[(\gamma_{n} + \gamma_{n+1} )(\Xi_n +\Xi_{n+1})\right] 
		\\&	\qquad \qquad 
		+\gamma_{n}\gamma_{n-1}\left(	\gamma_{n-2}+	\gamma_{n-1}+ \gamma_{n} + \gamma_{n+1}\right)\bigg]	\mathcal{S}_{n}(x;t) 
		\nonumber \\&\quad \nonumber
		+\gamma_{n}\gamma_{n-1}\bigg[ \Xi_{n-1}+ \Xi_n +\Xi_{n+1} +\gamma_{n-1}\gamma_{n+1}
		+\gamma_{n-2}\left( \gamma_{n-3}+\gamma_{n-2}+\gamma_{n-1}+  \gamma_{n} 
		+ \gamma_{n+1}\right)\bigg]
		\mathcal{S}_{n-2}(x;t) 
		\nonumber \\&\quad \nonumber
		+	\gamma_{n}  \gamma_{n-1} \gamma_{n-2} \gamma_{n-3} 	
		\bigg[
		\gamma_{n-4}+	\gamma_{n-3}+\gamma_{n-2}+\gamma_{n-1}+  \gamma_{n} + \gamma_{n+1}
		\bigg]\mathcal{S}_{n-4}(x;t)
		\nonumber \\&\quad 
		+	\left(\gamma_{n}  \gamma_{n-1} \gamma_{n-2} \gamma_{n-3}  \gamma_{n-4} \gamma_{n-5}\right)		\mathcal{S}_{n-6}(x;t),
		\end{align}}
	
	\begin{align}
	x^4 \mathcal{S}_n(x;t) &= \mathcal{S}_{n+4}(x;t) + \left(  \gamma_{n}+ \gamma_{n+1} + \gamma_{n+2} +\gamma_{n+3}\right)  \mathcal{S}_{n+2}(x;t)
	\nonumber\\& 
	+ \big[\gamma_{n} ( \gamma_{n-1} + \gamma_{n}  + \gamma_{n+1})    +\gamma_{n+1} ( \gamma_{n} + \gamma_{n+1}  + \gamma_{n+2} )\big]  \mathcal{S}_{n}(x;t)
	\nonumber\\& +\gamma_{n}  \gamma_{n-1}(  \gamma_{n-2} +\gamma_{n-1} + \gamma_{n}  + \gamma_{n+1} ) \mathcal{S}_{n-2}(x;t)
	+ ( \gamma_{n}  \gamma_{n-1} \gamma_{n-2} \gamma_{n-3} )  \mathcal{S}_{n-4}(x;t),
	\\
	x^2 \mathcal{S}_n(x;t) &= \mathcal{S}_{n+2}(x;t) + \left( \gamma_{n} + \gamma_{n+1}\right) \mathcal{S}_{n}(x;t) + \gamma_{n} \gamma_{n-1} \mathcal{S}_{n-2}(x;t),\label{recurrence1}
	\end{align}
\end{subequations}
\noindent
By substituting  Eq. \eqref{recurrence13a} into  Eq. \eqref{maineq}, we obtain the coefficients $\lbrace \mathfrak{f}_{n,j} \rbrace_{j=n-4}^{n-1}$ in  Eq.  \eqref{quasi} as:
{\small\begin{subequations}\label{Aau}
		\begin{align} 
		\mathfrak{u}_{n,n-6} &  =		6c	\left(\prod_{j=0}^{5}\gamma_{n-j}\right)= 6c\Big[\gamma_{n}  \gamma_{n-1} \gamma_{n-2}  \gamma_{n-3}\gamma_{n-4}  \gamma_{n-5}\Big], \,
		\qquad 
		\mathfrak{u}_{n,n-5}  =0,\\
		\mathfrak{u}_{n,n-4} &  =
		6c	\left(\prod_{j=0}^{3}\gamma_{n-j}\right)
		\Big[\gamma_{n-4} + \gamma_{n-3}+\gamma_{n-2} + \gamma_{n-1} +\gamma_{n} + \gamma_{n+1}\Big],
		\qquad 
		\mathfrak{u}_{n,n-3}  =0,
		\end{align}
		\begin{align}
		\nonumber\mathfrak{u}_{n,n-2} &= 
		\gamma_{n}  \gamma_{n-1}\Bigg[6c	\big\{ \Xi_{n-2} + \Xi_{n-1} +\Xi_{n} + \Xi_{n+1}
		+\gamma_{n-1} \gamma_{n-2} + \gamma_{n+1}\left(\gamma_{n-2} +  \gamma_{n-1}\right) \big\}
		\\&\qquad \qquad \qquad  +4t\left(\gamma_{n-2} + \gamma_{n-1} +\gamma_{n} + \gamma_{n+1}\right)-2t\Bigg]
		,\\
		\mathfrak{u}_{n,n-1} & =0.
		\end{align}
\end{subequations}}
For the case when $k = n$, we use integration technique    in  Eq. \eqref{Ccoeff} to obtain
{\small\begin{align}
	\Gamma_n \mathfrak{f}_{n,n}  &=\int_{-\infty}^{\infty}   x
	\dfrac{\,d\mathcal{S}_n (x;t)}{\,d{x}}
	\mathcal{S}_n(x;t)W_{\sigma}(x;t) \,\mathrm{d}x
	=-\tfrac12\int_{-\infty}^{\infty}   \mathcal{S}_n^2(x;t) \left[ W_{\sigma}(x;t) + x
	\dfrac{\,d W_{\sigma} (x;t)}{\,d{x}}
	\right] ~\mathrm{d} x
	\nonumber\\&= - \tfrac12\Gamma_n 
	+\int_{-\infty}^{\infty}   \mathcal{S}_n^2(x;t)
	\big( 3cx^{6}-2tx^{4}+tx^2
	-\sigma -\tfrac12\big)W_{\sigma}(x;t) \,\mathrm{d}x
	\nonumber\\&
	\nonumber=3c\int_{-\infty}^{\infty}  x^{6}\mathcal{S}_n^2(x;t)W_{\sigma}(x;t) \,\mathrm{d}x -2t\int_{-\infty}^{\infty}  x^{4}\mathcal{S}_n^2(x;t)W_{\sigma}(x;t) \,\mathrm{d}x
	\\& \qquad 
	+t\int_{-\infty}^{\infty}  x^2\mathcal{S}_n^2(x;t) W_{\sigma}(x;t) \,\mathrm{d}x  - (\sigma+1) \Gamma_{n}. 
	\label{eq:cnn1}
	\end{align}
	{\small\begin{subequations}
			\noindent
			By using the recursive relation given  in  Eq. \eqref{symmetricrecur} for  Eq. \eqref{eq:cnn1}, we  have that
			\begin{align}
			x^2 \mathcal{S}_n^2&=(\mathcal{S}_{n+1}+\gamma_n \mathcal{S}_{n-1})^2 
			= \mathcal{S}_{n+1}^2+2\gamma_n\mathcal{S}_{n+1}\mathcal{S}_{n-1}+\gamma_n^2\mathcal{S}_{n-1}^2,
			\\\nonumber
			x^4 \mathcal{S}_n^2&= x^2\big(\mathcal{S}_{n+1}^2+2\gamma_n\mathcal{S}_{n+1}\mathcal{S}_{n-1}+\gamma_n^2\mathcal{S}_{n-1}^2)\nonumber
			= x^2 \mathcal{S}_{n+1}^2 +2\gamma_n(x\mathcal{S}_{n+1})(x\mathcal{S}_{n-1})+\gamma_n^2x^2\mathcal{S}_{n-1}^2\\\nonumber
			&= \big(\mathcal{S}_{n+2}+\gamma_{n+1}\mathcal{S}_{n}\big)^2 + 2\gamma_n\big(\mathcal{S}_{n+2}+\gamma_{n+1}\mathcal{S}_{n}\big)\big(\mathcal{S}_{n}+\gamma_{n-1}\mathcal{S}_{n-2}\big)
			+\gamma_n^2\big(\mathcal{S}_{n}+\gamma_{n-1}\mathcal{S}_{n-2}\big) ^2\\ \nonumber
			&= \mathcal{S}_{n+2}^2+2(\gamma_{n+1}+\gamma_n)\mathcal{S}_{n+2}\mathcal{S}_{n}+(\gamma_{n+1}+\gamma_n)^2 \mathcal{S}_{n}^2+2\gamma_n\gamma_{n-1}\mathcal{S}_{n+2}\mathcal{S}_{n-2}
			\\ &\qquad\quad
			+2\gamma_n\gamma_{n-1}(\gamma_n+\gamma_{n+1})\mathcal{S}_{n}\mathcal{S}_{n-2}+\gamma_n^2\gamma_{n-1}^2\mathcal{S}_{n-2}^2,
			\end{align}
				\end{subequations}}
	and so by orthogonality, we have that
	{\small\begin{subequations}
			\begin{align}
			\nonumber
			\int_{-\infty}^{\infty} x^6 \mathcal{S}_n^2(x;t)~W_{\sigma}(x;t) \,\mathrm{d}x
			& =(\Gamma_{n+3}+\gamma_{n+2}^2\Gamma_{n+1}) +2(\gamma_{n+1}+\gamma_n)\gamma_{n+1}\gamma_{n+2}\Gamma_{n} 
			+(\gamma_{n+1}+\gamma_n)^2 (\Gamma_{n+1}+\gamma_{n}^2\Gamma_{n-1}) 
			\nonumber\\ & \qquad
			\nonumber\\ &= 
			(\gamma_{n}\gamma_{n+1}\gamma_{n+2})\Gamma_{n} + \Xi_{n+2}\gamma_{n+1}\Gamma_{n}
			+\gamma_{n+1}(\Xi_{n}+\Xi_{n+1})\Gamma_{n}+\gamma_{n}(\Xi_{n}+\Xi_{n+1})\Gamma_{n} 
			\nonumber\\ &\qquad\quad+ \gamma_{n-1} \gamma_n\gamma_{n+1}\Gamma_{n}
			+\gamma_{n}\Xi_{n-1} \Gamma_{n},
			\label{Neweq:cnn4}
			\end{align}
			
			\begin{align}
			\int_{-\infty}^{\infty}  x^2 ~\mathcal{S}_n^2(x;t)~W_{\sigma}(x;t) ~\,\mathrm{d}x \
			&= \Gamma_{n+1}+\gamma_n^2 ~\Gamma_{n-1} = (\gamma_{n+1}+\gamma_{n})~\Gamma_{n},\label{eq:cnn2} \\
			\int_{-\infty}^{\infty}  x^4 \mathcal{S}_n^2(x;t)~W_{\sigma}(x;t)~ \,\mathrm{d}x &= \Gamma_{n+2}+(\gamma_{n+1}+\gamma_n)^2\Gamma_n +\gamma_n^2\gamma_{n-1}^2~\Gamma_{n-2}
			\nonumber\\&=\big[(\gamma_{n+1}+\gamma_{n}+\gamma_{n-1})\gamma_n 
			+ (\gamma_{n+2}+\gamma_{n+1}+\gamma_n)\gamma_{n+1}\big]\Gamma_n\nonumber\\
			&= \left(\Xi_{n} +\Xi_{n+1}\right)\Gamma_n,
			\label{eq:cnn3}
			\end{align}
	\end{subequations}}
	using $\Gamma_{n+1}=\gamma_{n+1}\Gamma_n$, the difference equation  Eq. \eqref{nonlinerreucrequnforFreudAnewAx} and  $\Xi_n$ is given by  Eq. \eqref{greatqn}.
	
	\noindent
	By rearranging  Eq. \eqref{nonlinerreucrequnforFreudAnewAx} and taking $n \rightarrow n-1$ in  Eq. \eqref{nonlinerreucrequnforFreudAnewAx}, we have
	\begin{subequations}
			\begin{align}\label{nonlinarA}
			2t\Xi_{n}-t\gamma_{n}= \dfrac{n+ (2\sigma +1) \Omega_{n}}{2}
			-3c\left[
			\gamma_{n}\left( \Xi_{n-1} +\Xi_{n} +\Xi_{n+1}\right) + \gamma_{n-1} \gamma_{n}  \gamma_{n+1} \right],
			\end{align}
			\begin{align}\label{nonlinarAB}
			2t\Xi_{n+1}-t\gamma_{n+1}= \dfrac{n+1+ (2\sigma +1) \Omega_{n+1}}{2}
			-3c\left[
			\gamma_{n+1}\left( \Xi_{n} +\Xi_{n+1} +\Xi_{n+2}\right) + \gamma_{n} \gamma_{n+1}  \gamma_{n+2} \right]. 
			\end{align}
	\end{subequations}
	By combining  Eqs. \eqref{nonlinarA} and  Eq. \eqref{nonlinarAB}, we obtain 
	{\small\begin{align}
		\nonumber
		-2t\int_{-\infty}^{\infty} & x^4 \mathcal{S}_n^2(x;t)~W_{\sigma}(x;t) \,\mathrm{d}x +t\int_{-\infty}^{\infty}  x^2 \mathcal{S}_n^2(x;t)~W_{\sigma}(x;t) \,\mathrm{d}x 
		\nonumber\\ &=-\left(t\gamma_{n}-2t\Xi_{n}\right)-\left(t\gamma_{n+1}-2t\Xi_{n+1}\right)
		\nonumber\\ &= 
		-3c \Bigg[
		\gamma_{n}\left( \Xi_{n-1} +\Xi_{n} +\Xi_{n+1}\right) + \gamma_{n-1} \gamma_{n}  \gamma_{n+1}
		+\gamma_{n+1}\left( \Xi_{n} +\Xi_{n+1} +\Xi_{n+2}\right) + \gamma_{n} \gamma_{n+1}  \gamma_{n+2}\Bigg]
		\nonumber\\ &\qquad\quad  +  \dfrac{2n+1+ (2\sigma +1) \left(\Omega_{n}+\Omega_{n+1}\right)}{2}
		\nonumber\\ &= 
		-3c \Bigg[
		\gamma_{n}\left( \Xi_{n-1} +\Xi_{n} +\Xi_{n+1}\right) + \gamma_{n-1} \gamma_{n}  \gamma_{n+1}
		+\gamma_{n+1}\left( \Xi_{n} +\Xi_{n+1} +\Xi_{n+2}\right) + \gamma_{n} \gamma_{n+1}  \gamma_{n+2}\Bigg]
		\nonumber\\ &\qquad\quad  +  n+ (\sigma+1),
		\label{sweet}
		\end{align}}
	Hence from  Eq. \eqref{Neweq:cnn4} and  Eq. \eqref{sweet},  Eq. \eqref{eq:cnn1} becomes
	{\small\begin{align} 
		\nonumber
		\mathfrak{u}_{n,n}&=
		\dfrac{1}{\Gamma_n} \biggl\{
		3c\int_{-\infty}^{\infty}  x^{6}\mathcal{S}_n^2(x;t)\,W_{\sigma}(x;t) \,\mathrm{d}x - (\sigma+1) \Gamma_{n} 
		-2t\int_{-\infty}^{\infty}  x^{4}\mathcal{S}_n^2(x;t)\,W_{\sigma}(x;t) \,\mathrm{d}x
		\nonumber\\ &\qquad	\qquad	
		+t\int_{-\infty}^{\infty}  x^2\mathcal{S}_n^2(x;t)\,W_{\sigma}(x;t) \,\mathrm{d}x 
		\biggr\}
		\nonumber\\ &= 
		3c \Bigg[ (\gamma_{n}\gamma_{n+1}\gamma_{n+2})+ \Xi_{n+2}\gamma_{n+1}
		+\gamma_{n+1}(\Xi_{n}+\Xi_{n+1})\Gamma_{n}+\gamma_{n}(\Xi_{n}+\Xi_{n+1}) 
		+ \gamma_{n-1} \gamma_n\gamma_{n+1}\Gamma_{n}+\gamma_{n}\Xi_{n-1} \Bigg] 
		\nonumber\\ & - (\sigma+1) 
		-3c \Bigg[ \gamma_{n}\left( \Xi_{n-1} +\Xi_{n} +\Xi_{n+1}\right) + \gamma_{n-1} \gamma_{n}  \gamma_{n+1}
		+\gamma_{n+1}\left( \Xi_{n} +\Xi_{n+1} +\Xi_{n+2}\right) + \gamma_{n} \gamma_{n+1}  \gamma_{n+2}\Bigg]
		\nonumber\\ &\quad  +  n+ (\sigma+1)
		\nonumber\\ & 
		=n. \label{eq:cnn4}
		\end{align}}
	Combining  Eq. \eqref{Aau} with  Eq. \eqref{quasi} gives
	{\small\begin{subequations}
			\begin{align}\label{AQquasi}
			x\dfrac{\,d\mathcal{S}_n}{~\mathrm{d} x}
			=\mathfrak{u}_{n,n-6} ~\mathcal{S}_{n-6}(x;t) +\mathfrak{u}_{n,n-4}~ \mathcal{S}_{n-4}(x;t) +  \mathfrak{u}_{n,n-2} ~\mathcal{S}_{n-2}(x;t)  +  \mathfrak{u}_{n,n}~ \mathcal{S}_{n}(x;t).
			\end{align}
			\noindent  Rewriting  $\mathcal{S}_{n-4}$ and $\mathcal{S}_{n-2}$ into Eq. \eqref{AQquasi} in terms of $\mathcal{S}_n$ and $\mathcal{S}_{n-1}$ using Eq. \eqref{symmetricrecur}, we obtain
			\begin{align}\label{recoa}
			\mathcal{S}_{n-2}(x;t)&= \dfrac{x\mathcal{S}_{n-1}(x;t)-\mathcal{S}_n(x;t)}{\gamma_{n-1}}, \\
			\mathcal{S}_{n-3}(x;t)&= \dfrac{x\mathcal{S}_{n-2}(x;t)-\mathcal{S}_{n-1}(x;t)}{\gamma_{n-2}} 
			= \dfrac{x^2-\gamma_{n-1}}{\gamma_{n-1}\gamma_{n-2}}\, \mathcal{S}_{n-1}(x;t) - \dfrac{x}{\gamma_{n-1}\gamma_{n-2}}\, \mathcal{S}_{n}(x;t),\\ 
			\mathcal{S}_{n-4}(x;t)
			&= \dfrac{x\mathcal{S}_{n-3}(x;t) - \mathcal{S}_{n-2}(x;t)}{\gamma_{n-3}} 
			= \label{recoab}   \dfrac{x^3-  ( \gamma_{n-1}+\gamma_{n-2})x}{\gamma_{n-1}\gamma_{n-2}\gamma_{n-3}}\,  \mathcal{S}_{n-1}(x;t) -
			\dfrac{x^2-\gamma_{n-2}}{ \gamma_{n-1}\gamma_{n-2}\gamma_{n-3}} \,  \mathcal{S}_{n}(x;t),
			\\
			\mathcal{S}_{n-6}(x;t)
			\nonumber 
			= &
			\Biggl\{
			\dfrac{x^5-( \gamma_{n-1}+\gamma_{n-2}+\gamma_{n-3}+\gamma_{n-4})x +\left(\gamma_{n-1}\gamma_{n-3}+\gamma_{n-1}\gamma_{n-4}+\gamma_{n-2}\gamma_{n-4}\right)}{\gamma_{n-1}\gamma_{n-2}\gamma_{n-3}\gamma_{n-4}
				\gamma_{n-5}
			} 
			\Biggr\}
			\mathcal{S}_{n-1}(x;t)
			\nonumber\\ &\qquad-\Biggl\{
			\dfrac{x^4-(\gamma_{n-2}+\gamma_{n-3}+\gamma_{n-4})x +\gamma_{n-2}\gamma_{n-4}}{\gamma_{n-1}\gamma_{n-2}\gamma_{n-3}\gamma_{n-4}\gamma_{n-5}} 
			\Biggr\}
			\mathcal{S}_{n}(x;t).
			\label{Hexaky}
			\end{align}
	\end{subequations}}
	\noindent Substituting  Eq. \eqref{Aau},   Eq. \eqref{eq:cnn4},  Eq. \eqref{recoa}, Eq. \eqref{recoab} and Eq. \eqref{Hexaky} into  Eq. \eqref{AQquasi} yields  the required result.
	\section{The differential equation}
	\label{secdiffFreud}
	\noindent
	\begin{theorem}\label{DiffeqnSextic}
		For the semiclassical weight in \eqref{sexticFreud}, the corresponding monic orthogonal polynomials $\mathcal{S}_{n}(x;t)$ obey a linear ODE (with rational coefficients) as
		\begin{equation}
		\label{eq:SnodeFreud}
		\dfrac{\,d^2}{\,d{x^2}}
		\mathcal{S}_{n}(x;t)+\tilde{U}_n(x;t)~
		\dfrac{\,d}{\mathrm{d}x}
		\mathcal{S}_{n}(x;t)
		+\tilde{W}_n(x;t)~\mathcal{S}_{n}(x;t)=0,
		\end{equation}
		where
		{\small\begin{subequations}\label{hexicdiffcoefffwwa}
				\begin{align}
				\nonumber
				\tilde{U}_n(x;t) 
				&=	-6cx^5-t(4x^3-2x)+\dfrac{(2\sigma +1)}{x}
				\nonumber\\ &\qquad
				-\Bigg[\dfrac{24cx^3+2\left[6c(\gamma_n+\gamma_{n+1}) +4t\right]x}{6cx^4+6c(\gamma_n+\gamma_{n+1})x^2+ 6c\left(\Xi_{n+1}+\Xi_{n}\right)-2t+4t\left(x^2+\gamma_n +\gamma_{n+1}\right)}\Bigg] 
				\end{align}
				{\footnotesize\begin{align}
					\nonumber	
					\tilde{W}_n(x;t) &= 18c\gamma_nx^2 + 6c\Xi_n -\dfrac{(2\sigma +1)\Omega_{n}}{x^2} +4t\gamma_n
					\nonumber\\&  \qquad 
					+\gamma_n\Bigg(
					6cx^4+6c(\gamma_n+\gamma_{n+1})x^2+ 6c\left(\Xi_{n+1}+\Xi_{n}\right)-2t+4t\left(x^2+\gamma_n +\gamma_{n+1}\right)
					\Bigg)
					\nonumber \\& \qquad \qquad	\qquad
					\times
					\Bigg(
					6cx^4+6c(\gamma_n+\gamma_{n-1})x^2+ 6c\left(\Xi_{n-1}+\Xi_{n}\right)-2t+4t\left(x^2+\gamma_n +\gamma_{n-1}\right)	\Bigg)
					\nonumber \\& \qquad  
					-\Bigg[
					\Big(6cx^5+(6c\gamma_n+4t)x^3- \dfrac{2\sigma +1}{x}+ \left(6c\Xi_n +4t\gamma_n -2t\right)x  +\dfrac{(2\sigma +1)\Omega_{n}}{x} 
					\nonumber \\& \qquad \qquad  \qquad 
					+\dfrac{24cx^3+2\left[6c(\gamma_n+\gamma_{n+1}) +4t\right]x}{6cx^4+6c(\gamma_n+\gamma_{n+1})x^2+ 6c\left(\Xi_{n+1}+\Xi_{n}\right)-2t+4t\left(x^2+\gamma_n +\gamma_{n+1}\right)}
					\Big)
					\nonumber \\& \qquad   \qquad 
					\times \Big(6c\gamma_nx^3 +\left(6c\Xi_n +4t\gamma_n\right)x +\dfrac{(2\sigma +1)\Omega_{n}}{x}\Big)\Bigg]
					\nonumber
					\\& \equiv
					- \mathcal{B}_n(x;t)\left[v'(x) +\mathcal{B}_n(x;t) +\frac{\mathcal{A}'_n (x;t)}{ \mathcal{A}_n(x;t)}\right] + \gamma_n \mathcal{A}_n(x;t) A_{n - 1}(x;t)+	\mathcal{B}'_n (x;t),		
					\end{align}}
		\end{subequations}}
		where {\small$\Omega_{n}$} and  {\small$\Xi_n$}  are given in Eqs.  \eqref{OmegaN} and \eqref{greatqn}  respectively. 
	\end{theorem}
	\begin{proof}
		For  the proof, consult similar ideas in  \cite{kelil2018properties} and \cite{refkelilAppadu2020}.
	\end{proof}

\begin{remark}
	One can expand  Eq. \eqref{hexicdiffcoefffwwa} via symbolic packages such as Mathematica (Maple), however the resulting expression  may look quite cumbersome.  
\end{remark}

\section{Application of Eq. \eqref{eq:SnodeFreud} for electrostatic zero distribution}
\noindent The authors in \cite{garrido2005electrostatic} considered  a perturbation of quartic Freud weight ($w(x)=\exp(-x^4)$) by the addition of a fixed charged point of mass {\small$\delta$} at the origin; the corresponding polynomials are Freud-type polynomials (see the recent work in \cite{refGarzaHuertasMarcellan}). For semiclassical orthogonality measure, it was shown in \cite{garrido2005electrostatic} that these polynomials obey a second-order linear differential equation of the form \eqref{H}, and the electrostatic model is in sight as in  \cite{refIsmaila}. Application of Eq. \eqref{eq:SnodeFreud} for electrostatic zero distribution is also mentioned. 
\noindent Following these ideas, a similar work for the perturbed Freud-type weight in \eqref{sexticFreud} is given  in a recent paper \cite{refkelilAppadu2020} using  the obtained differential equation in Section  \ref{secdiffFreud}.

\section{Conclusions}
\noindent
By introducing a time variable to scaled sextic Freud-type measure upon deformation (perturbation), we have found certain fresh characterizing properties: some recursive relations, moments of finite order, concise formulation and orthogonality relation, nonlinear difference equation for recurrence coefficients as well as the corresponding polynomials and certain properties of the zeros of the corresponding polynomials. This work derived certain non-linear difference equations, Toda-like equations, and differential equations for the recurrence coefficients of the corresponding orthogonal polynomials under consideration. 
Special attention, using the method of Shohat's quasi-orthogonality and ladder operators, is given 
to characterize the Freud-type weight \eqref{sexticFreud}. Such semiclassical symmetric weight in \eqref{sexticFreud} follows from quadratic transformation and symmetrization  as in \cite{refChihara}.
By combining  the three-term recurrence relation with the difference-recurrence relation, a second-order differential equation fulfilled by polynomials associated with the semiclassical weight \eqref{sexticFreud} is obtained. Application of the resulting differential equation in Eq. \eqref{eq:SnodeFreud} for electrostatic zero distribution is also noted.  Following this work, investigation of these recurrence coefficients in  connection with certain (discrete) integrable systems will be a prominent continuation of this  study.
%
%

\end{document}